\documentclass[a4paper,leqno]{siamltex}

\usepackage[thinspace,amssymb,thinqspace,thinspace]{SIunits}
\usepackage{amsmath}
\usepackage{stmaryrd}
\usepackage{amssymb}
\usepackage{graphicx}
\usepackage{subfigure}
\usepackage{xspace}
\usepackage[noend]{algorithmic}
\usepackage{cite}
\usepackage{soul}

\usepackage[ruled]{algorithm2e}

\usepackage[left=3cm,right=2cm,top=2cm,bottom=2cm]{geometry}
\usepackage[utf8]{inputenc}

\usepackage{tikz}

\usepackage[citecolor=black]{hyperref}

\newcommand{\VectDiscr}[1]{\,{#1}}
\newcommand{\MatrDiscr}[1]{{#1}}
\newcommand{\VD}[1]{\VectDiscr{#1}}
\newcommand{\MD}[1]{\MatrDiscr{#1}}

\newcommand{\muelu}{\textsc{MueLu}}
\newcommand{\ifpackTwo}{\textsc{Ifpack2}}
\newcommand{\zoltanTwo}{\textsc{Zoltan2}}

\newcommand{\todo}[2]{\textcolor{red}{{\normalfont\bfseries ToDo (#1):} #2}}

\newcommand{\verify}[1]{\textcolor{magenta}{{\normalfont\bfseries Verify:} #1}}

\newcommand{\REMOVE}[1]{}

\newcommand{\rstumin}[1]{\textcolor{olive}{Ray: #1}}
\newcommand{\mmayr}[1]{\textcolor{cyan}{Matthias: #1}}

\newcommand{\secref}[1]{Section~\ref{#1}} 
\newcommand{\Secref}[1]{Section~\ref{#1}} 
\newcommand{\Tabref}[1]{Table~\ref{#1}} 
\newcommand{\algref}[1]{Algorithm~\ref{#1}} 

\newtheorem{remark}{Remark}[section]

\newcommand{\define}[1]{\emph{#1}} 

\newcommand{\multigrid}{multigrid}
\newcommand{\Multigrid}{Multigrid}

\newcommand{\bbox}{black--box} 
\newcommand{\vcycle}{V-cycle}
\newcommand{\wcycle}{W-cycle}

\newcommand{\vs}{vs.}
\newcommand{\wrt}{w.r.t.}

\newcommand{\GS}{Gauss--Seidel} 
\newcommand{\Cheby}{Chebyshev} 

\newcommand{\dom}{\Omega} 

\newcommand{\AnyQuantity}{\left(\bullet\right)}

\newcommand{\sub}[2]{#1_{#2}} 

\newcommand{\indLevel}{\ell} 

\newcommand{\inv}[1]{#1^{-1}} 
\newcommand{\trans}[1]{#1^{\mathrm{T}}} 

\newcommand{\matrixBlock}[2]{(#1,#2)} 

\newcommand{\restrictor}{R} 
\newcommand{\prolongator}{P} 
\newcommand{\RAP}{\Coarse{\linMat}} 
\newcommand{\smoother}{\mathcal{S}} 
\newcommand{\indexedLevel}[2]{#1_{#2}} 
\newcommand{\levelL}[1]{\indexedLevel{#1}{\indLevel}} 
\newcommand{\coarse}{\mathrm{c}}

\newcommand{\Coarse}[1]{{#1}_{\coarse}}
\newcommand{\numLevels}{L} 
\newcommand{\numLevelsRegion}{L_{\hhgRegional}} 
\newcommand{\numLevelsAMG}{L_{\hhgComposite}} 
\newcommand{\numPreSmooth}{\nu_1} 
\newcommand{\error}{e}

\newcommand{\linMat}{A} 
\newcommand{\linSol}{x} 
\newcommand{\linRhs}{b} 
\newcommand{\linVec}{v} 

\newcommand{\boostFactor}{\kappa}
\newcommand{\tildeEigValMax}{\tilde{\lambda}_{\mathrm{max}}} 
\newcommand{\eigValMax}{\lambda_{\mathrm{max}}} 
\newcommand{\eigValMin}{\lambda_{\mathrm{min}}} 
\newcommand{\eigRatio}{\eta} 

\newcommand{\region}[2]{{#1}^{(#2)}} 
\newcommand{\regionSub}[3]{{#1}_{#3}^{(#2)}} 
\newcommand{\hhgComposite}{\mathrm{c}}
\newcommand{\HHGComposite}[1]{{#1}}
\newcommand{\hhgRegional}{\mathrm{r}}
\newcommand{\HHGRepRegional}[1]{\llbracket {#1} \rrbracket }
\newcommand{\HHGRegional}[1] {\hskip .01in {\bf [ \hskip -.045in [ \hskip -.045in [} \hskip .01in {#1} \hskip .01in {\bf ] \hskip -.045in ] \hskip -.045in ]}\hskip .01in}
\newcommand{\indReg}{k} 

\newcommand{\Int}{\mathrm{I}} 
\newcommand{\Bound}{\mathrm{\Gamma}} 
\newcommand{\Bounds}[1]{\Bound_{#1}} 
\newcommand{\Boundss}[2]{\Bound_{#1#2}} 

\newcommand{\regionIntOp}{\mathcal{R}} 
\newcommand{\RegionIntOp}[1]{\regionIntOp\left(#1\right)} 
\newcommand{\regToComp}{\mathcal{C}} 
\newcommand{\RegToComp}[1]{\regToComp\left(#1\right)} 
\newcommand{\cPsi}{\bar{\Psi}}
\newcommand{\nRegions}{m}

\newcommand{\nproc}{n^{\mathrm{proc}}} 

\usetikzlibrary{arrows.meta}
\usetikzlibrary{bending}
\usetikzlibrary{positioning}
\usetikzlibrary{snakes}
\usetikzlibrary{shapes.misc}
\usetikzlibrary{shapes.geometric}

\tikzstyle{domainline}=[very thick]
\tikzstyle{hiddenline}=[thick,dashed]
\tikzstyle{dbcline}=[thick]
\tikzstyle{tractionline}=[thick]
\tikzstyle{tractionarrow}=[thick,-{Latex}]
\tikzstyle{measline}=[thick, {Latex}-{Latex}]
\tikzstyle{measradius}=[thick, -{Latex}]
\tikzstyle{measauxline}=[thick]
\tikzstyle{coordsline}=[thick,{Latex[open,fill=white]}-{Latex[open,fill=white]}]
\tikzstyle{coordline}=[thick,-{Latex[open,fill=white]}]
\tikzstyle{symmetryline}=[very thick,dash pattern=on 10pt off 3pt on \the\pgflinewidth off 3pt]

\tikzstyle{darkgreen}=[green!50!black]

\tikzstyle{springlineal}=[snake=zigzag,thick,line before snake=0.5cm,line after snake=0.5cm,segment length=6,segment amplitude=5,join=round]

\tikzstyle{flowchartarrow}=[thick,-{Latex}]

\graphicspath{{fig/}{eps/}}

\title{Non-invasive multigrid for semi-structured grids\thanks{
This work was supported by the U.S.~Department of Energy, Office of Science, Office of Advanced Scientific Computing Research, Applied Mathematics program.  Sandia National Laboratories is a multimission laboratory managed and operated by National Technology and Engineering Solutions of Sandia, LLC., a wholly owned subsidiary of Honeywell International, Inc., for the U.S. Department of Energy's National Nuclear Security Administration under grant DE-NA-0003525. This paper describes objective technical results and analysis. Any subjective views or opinions that might be expressed in the paper do not necessarily represent the views of the U.S. Department of Energy or the United States Government. SAND2021-3211 O
}}
\author{Matthias Mayr\thanks{Institute for Mathematics and Computer-Based Simulation, University of the Bundeswehr Munich, Werner-Heisenberg-Weg 39, 85577 Neubiberg, Germany (matthias.mayr@unibw.de), \emph{This work was partially performed while this authors was affiliated with Sandia National Laboratories, Livermore, CA 94551}},
        \and 
        Luc Berger-Vergiat\thanks{Sandia National Laboratories, Albuquerque, NM 87185 (lberge@sandia.gov)},
        \and 
        Peter Ohm\thanks{Sandia National Laboratories, Albuquerque, NM 87185 (pohm@sandia.gov)}, 
        \and 
        Raymond S. Tuminaro\thanks{Sandia National Laboratories, Livermore, CA 94551 (rstumin@sandia.gov)}}
\date{Today} 

\begin{document}
\maketitle

\begin{abstract}
Multigrid solvers for hierarchical hybrid grids (HHG) have been proposed 
to promote the efficient utilization of high performance 
computer architectures.  These HHG meshes are constructed by 
uniformly refining a relatively coarse fully unstructured mesh.  
While HHG meshes provide some flexibility for unstructured applications, 
most multigrid calculations can be accomplished using efficient structured
grid ideas and kernels.  This paper focuses on generalizing the HHG idea 
so that it is applicable to a broader community of computational scientists,
and so that it is easier for existing applications to leverage structured 
multigrid components. Specifically, we adapt the structured multigrid methodology
to significantly more complex semi-structured meshes. Further, we illustrate
how mature applications might adopt a semi-structured solver in
a relatively non-invasive fashion. To do this, we propose a formal mathematical
framework for describing the semi-structured solver. This formalism allows
us to precisely define the associated multigrid method and to show its
relationship to a more traditional multigrid solver. Additionally, the 
mathematical framework clarifies the associated software design and implementation.
Numerical experiments highlight the relationship of the new solver with classical multigrid.
We also demonstrate the generality and potential performance gains associated
with this type of semi-structured multigrid. 
\end{abstract}

\section{Introduction}
\label{sec:Introduction}

\begin{figure}[ht!]
  \centering
  \includegraphics[scale=0.1]{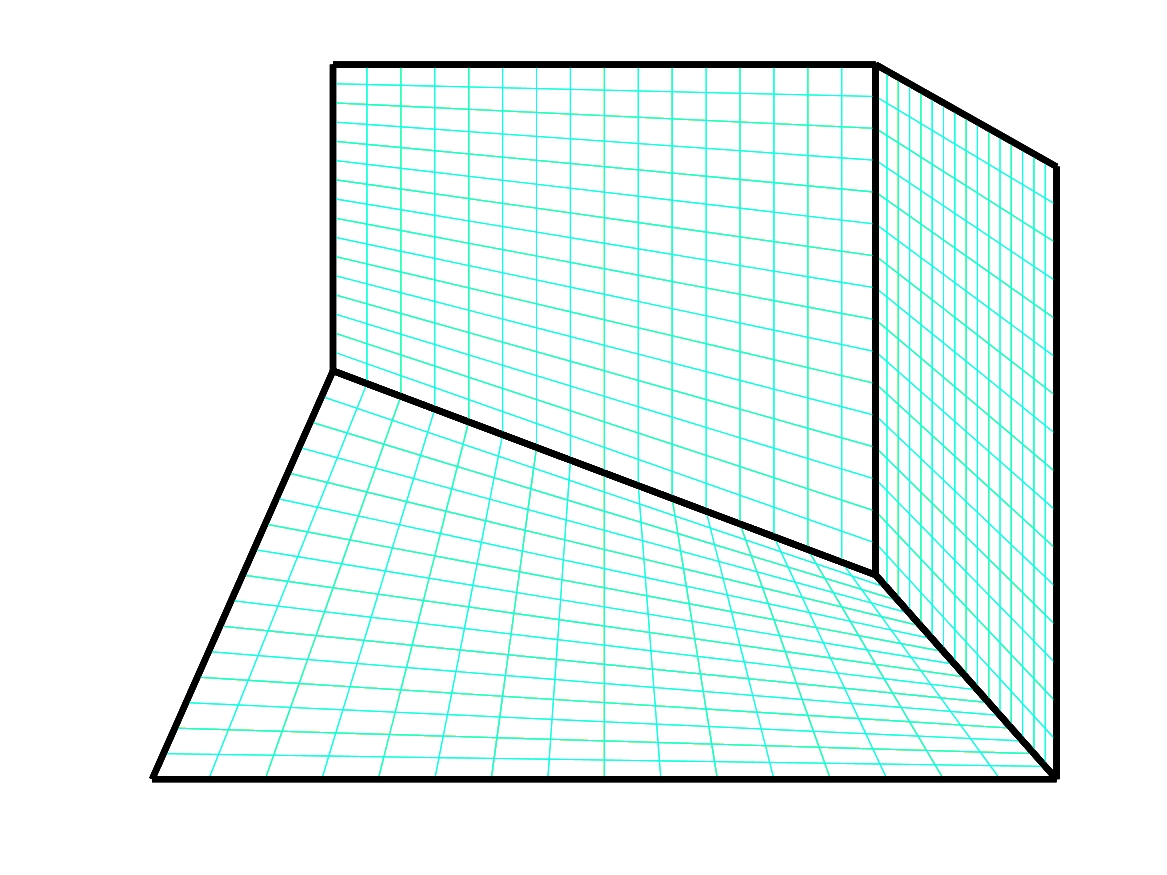}
  \includegraphics[scale=0.1]{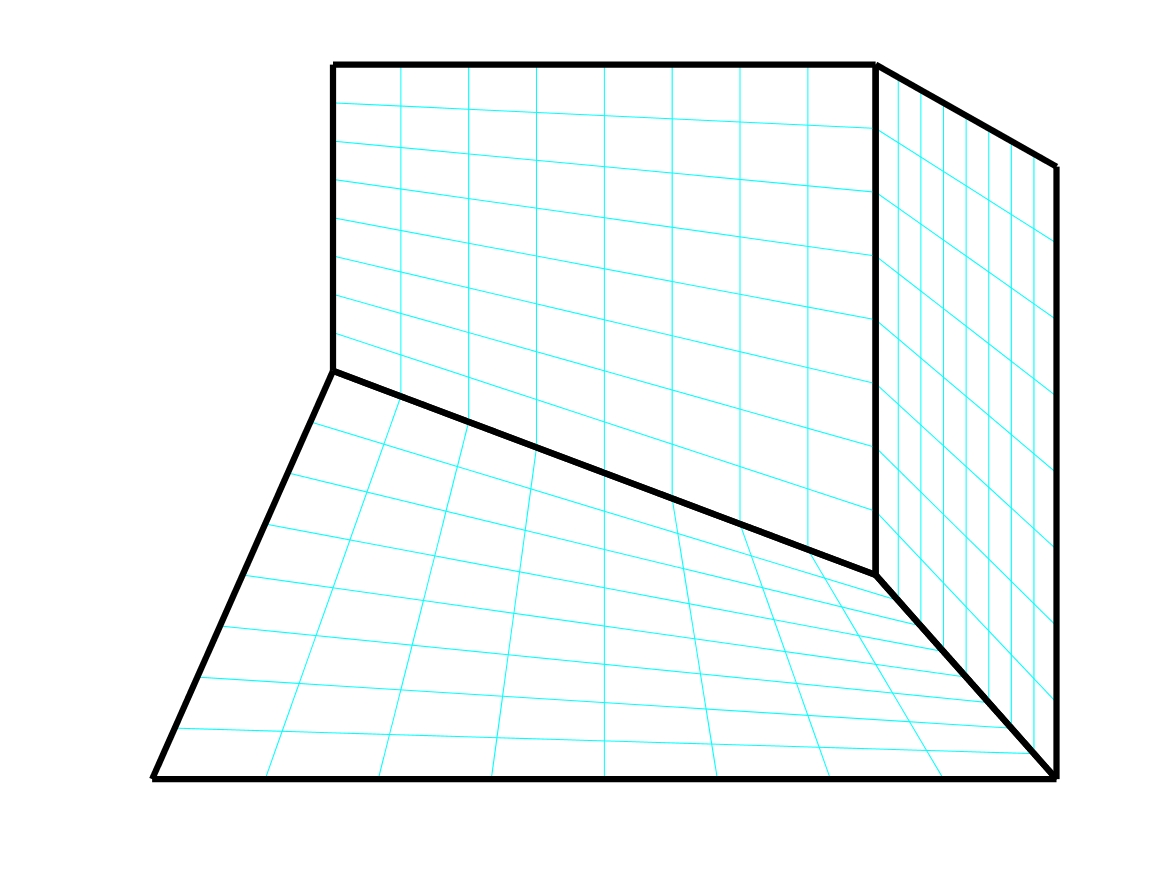}
  \includegraphics[scale=0.1]{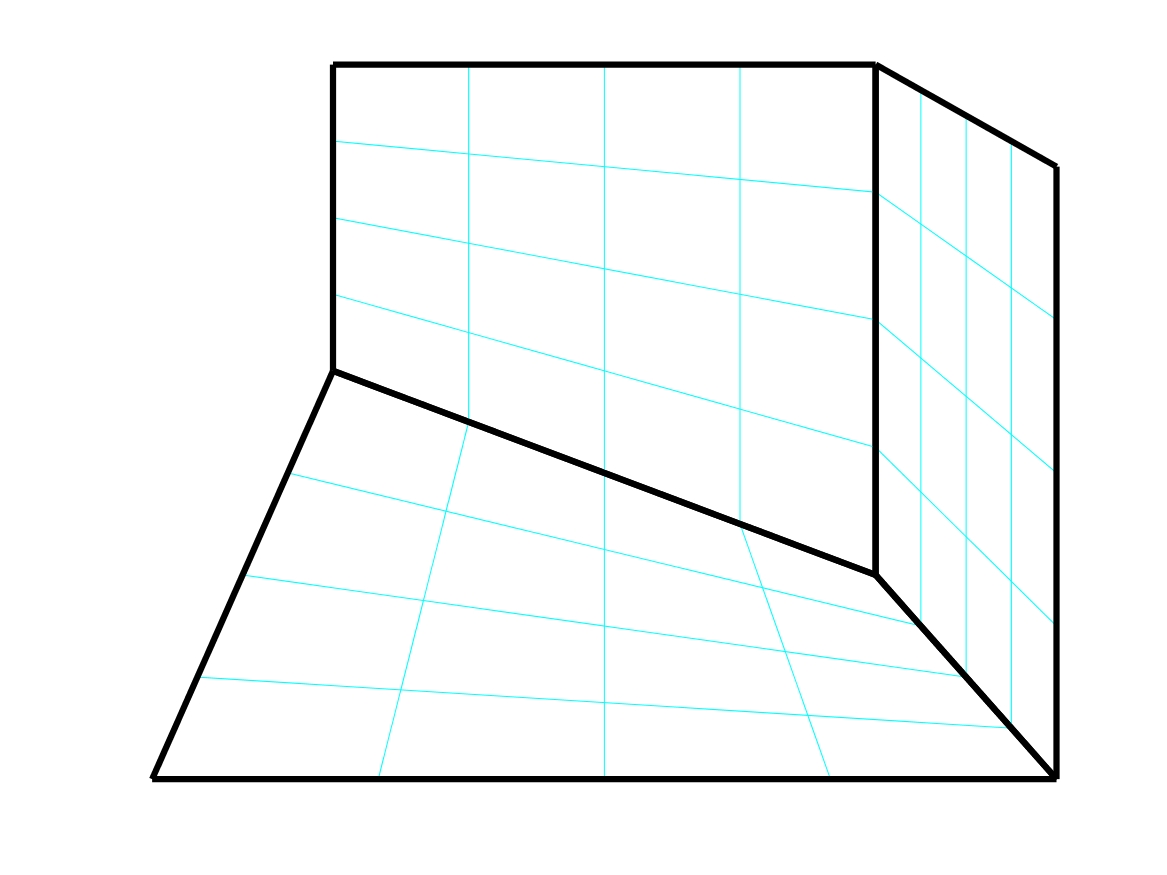}
  \caption{A hierarchy of two dimensional HHG meshes created by regular refinement of a 3 element mesh}
  \label{PBO:fig:meshRefine}
\end{figure}

Multigrid (MG) methods have been developed for both structured and unstructured
grids~\cite{Briggs2000a,Hackbusch1994a,Saad2003a,TrOoSc00}. In general, unstructured meshes are heavily favored
within sophisticated science and engineering simulations as they facilitate
the representation of complex geometric features. While unstructured 
approaches are often convenient, there are significant
potential advantages to structured meshes on exascale systems
in terms of memory, setup time, and kernel optimization.
In recent years, multigrid solvers for hierarchical hybrid grids (HHGs) have 
been 
proposed to provide some flexibility for unstructured applications 
while also leveraging some features of structured multigrid for performance
on advanced computing systems~\cite{Bergen2004a}. Hierarchical hybrid grids
are formed by regular refinement of an initial coarse grid.  The result 
is a HHG grid hierarchy containing regions of structured mesh, even if the 
initial coarse mesh is completely unstructured~\cite{Bergen2004a}. Essentially,
each structured region in an HHG mesh corresponds to one element of the 
original coarse mesh that has been uniformly refined.
A corresponding multigrid solver can then be developed using primarily structured
multigrid ideas.  Figure~\ref{PBO:fig:meshRefine} illustrates a two 
dimensional HHG mesh hierarchy with three structured 
regions.
Here,  the two rightmost grids might be used as multigrid
coarse grids for a discretization on the finest mesh.
The key point is that structured multigrid kernels can be used for most of the computation. These structured computations require significantly less memory and generally less communication than their unstructured counterparts. Further, the structured multigrid kernels are significantly 
more amenable to performance optimization on advanced architectures.
A series of papers~\cite{Bergen2006a,Bergen2007a,Gmeiner2012a,Gmeiner2013a,Gmeiner2015b,Gmeiner2015c} have documented
noticeably impressive HPC performance using an HHG approach on realistic
simulations, some involving over one trillion unknowns. In these papers, 
the primarily structured nature of the mesh is heavily leveraged throughout
the multigrid solver in an essentially matrix-free fashion.

While HHG solvers provide some balance between flexibility and structured performance, 
they do impose restrictions on the type of meshes that can be considered. 
Additionally, it is difficult to adapt existing finite element applications
to HHG solvers.
Of course there are alternative approaches to structure including 
composite grids, overset meshes, and octree meshes 
(for example~\cite{Philip20122277,Dubey20143217,Henshaw20087469,octrees,Lee03asynchronousfast,uintah}).
Additionally, \textsc{Hypre} has some semi-structured capabilities~\cite{hypre}.
While these approaches can also attain good scalability on high performance architectures, 
most scientific teams have been resistant to investigate these
structured grid possibilities
due to concerns about their intrusive nature, often requiring fundamental changes to 
the mesh representations and discretization technology employed within the application.  
This is especially true for unstructured finite element simulations, which dominate 
the discretization approaches employed at Sandia.

Our aim in this paper is to at least partially address these obstacles
by broadening the HHG approach to a wider class of meshes 
and by providing an easier or less-invasive code path to migrate
existing applications toward semi-structured solvers. To do this, we introduce
a mathematical framework centered around the idea of a {\it region} representation.
The region perspective decomposes the original domain into a set of regions that 
only overlap at inter-region interfaces and where the computational mesh also conforms
at these interfaces. The main difference from the typical situation
(which we refer to as the {\it composite} mesh to emphasize the differences)
is that each region has its own copy of solution unknowns along its interfaces. 
If all regions are structured, the overall 
grid is a block structured mesh (BSM). BSMs can be constructed by joining
separately meshed components or a regular refinement of an unstructured mesh as in
the HHG case. Thus, BSMs are a generalization of the HHG idea.
As in the HHG case, a special region-oriented solver can take advantage of structure
within structured regions. 

The mathematical framework allows us to consider region-oriented versions of algorithms 
developed from a traditional composite mesh perspective. It also provides conditions
on the region-oriented grid transfer operators to guarantee a  mathematical equivalence 
relationship between region-oriented multigrid and a traditional solver. In some cases, it is
easy to accomplish this exact equivalence while in other cases there are practical tradeoffs
that must be weighed, comparing additional computational/communication requirements
against a possible convergence benefit to exact equivalence. One key result of the mathematical 
framework
is that in some cases (linear interpolation grid transfers without curved region interfaces) it is 
possible to construct a region multigrid hierarchy without communication.
This includes no communication requirement for the Galerkin triple matrix product 
(used to project the discretization operator) when all associated matrices adopt a region
representation. This is in contrast to a standard AMG setup algorithm where communication
costs can be noticeable especially when the density of the discretization sparsity pattern increases
as one constructs coarser and coarser matrices. 

\begin{figure}[ht!]
\centering
\hfill
\includegraphics[trim = 5.0in 1.3in 0.6in 2.6in, clip, height = 3.1cm,width = 3.5cm]{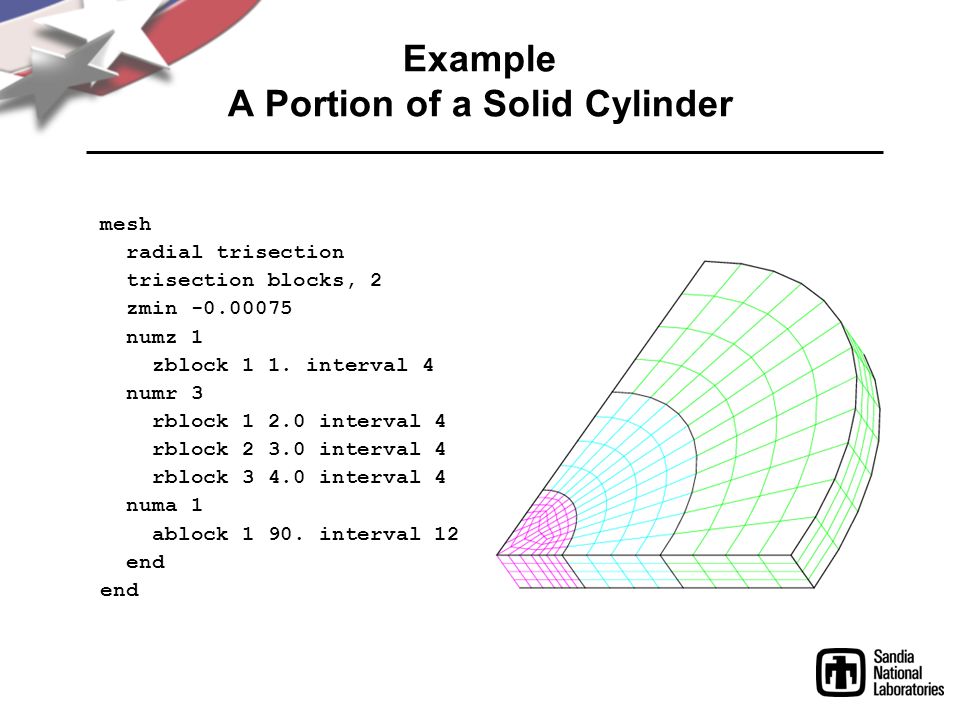}
\hfill
\includegraphics[trim = 3.5in 2.3in 3.2in 3.0in, clip, height = 3.1cm,width = 3.5cm]{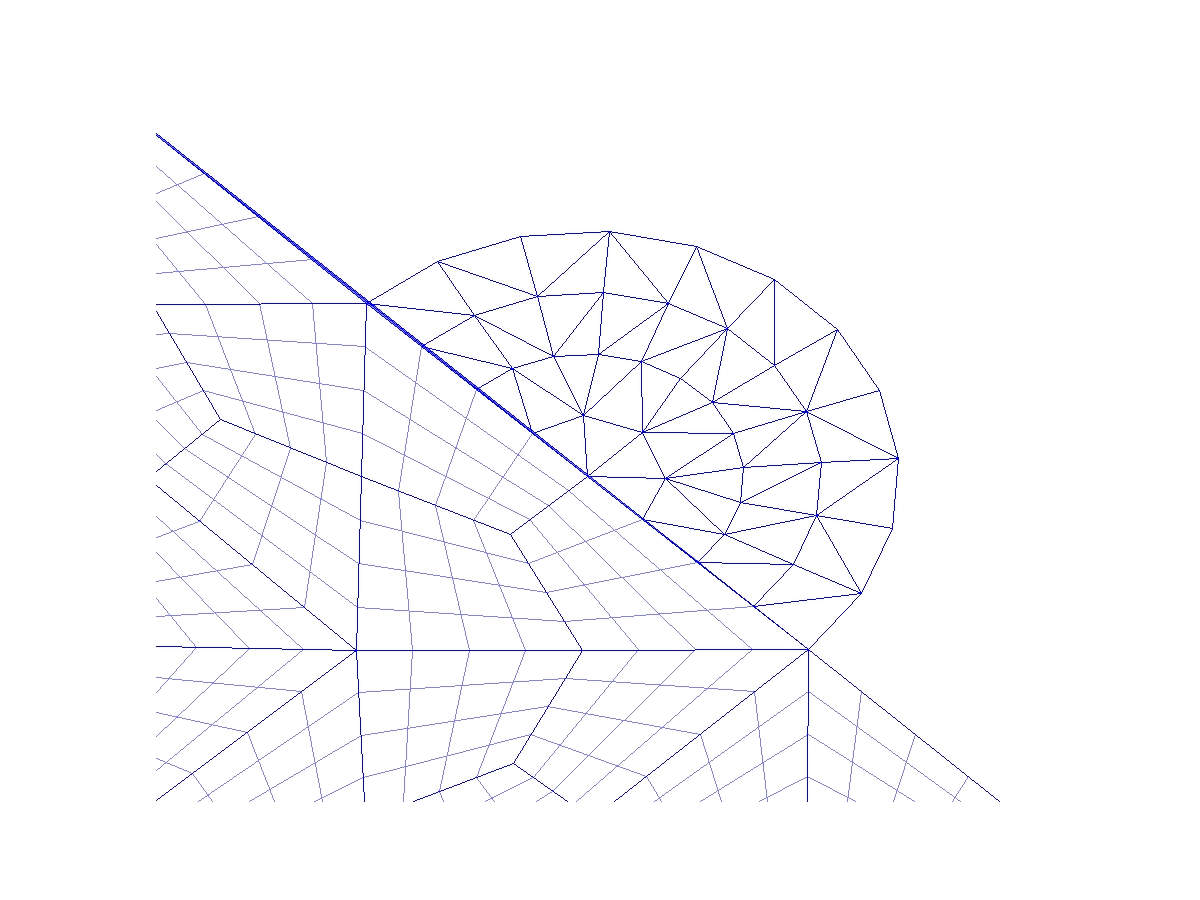}
\hfill
\includegraphics[trim = 1.2in 0.0in 0.3in 0.1in, clip, height = 3.9cm,width = 3cm]{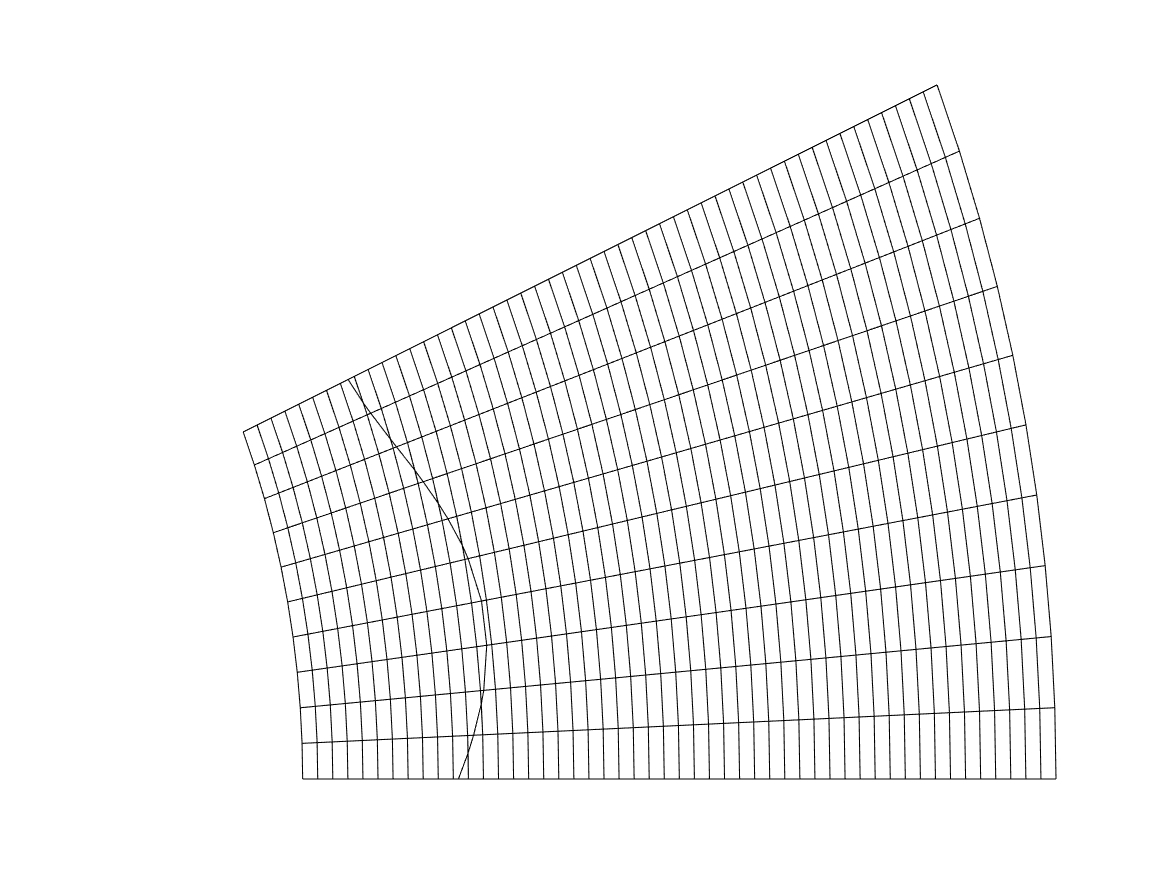}
\hfill~
\caption{Radial tri-section mesh (left), unstructured region attached to an HHG mesh (middle), interface with cut element mesh (right).\label{radtri}}
\end{figure}

The mathematical framework is fairly general in that it is not restricted to structured regions. 
That is, it allows for the possibility that some regions might be structured while
others are unstructured. This can be useful in applications where it might be awkward 
to  resolve certain geometries or to capture local features with only structured regions.
Figure~\ref{radtri}
illustrates some partially structured meshes.  The leftmost image corresponds to a mesh
used to represent wires.  The middle picture illustrates a main body mesh with an attached part.
The rightmost example displays a background mesh with some split elements to represent
an interface. In this last case, an unstructured region might be employed only to
surround the interface. 
Our software considers these types of situation again
using the mathematical framework as a guide for the treatment of 
grid transfer operators near region interfaces.
Of course, a  matrix-free approach would be problematic in this more general setting 
and performance in unstructured regions might be poorer, though there will be much fewer 
unstructured regions. 

One nice aspect of the mathematical framework is that it formalizes the transformation
between composite and region perspectives. As noted, this is helpful when designing
grid transfers near region interfaces. It is also helpful, however, when understanding
the minimal application requirements for employing such a region-oriented solver.
In particular, the finite element software must provide a structured PDE 
matrix for each structured region as well as more detailed information on how to glue
regions together.  It is easy for the requirements of a semi-structured or an HHG framework
to become intrusive on the application infrastructure.
The philosophy taken in this paper is toward the development
of algorithms and abstractions that are sufficiently flexible to model 
complex features without imposing over-burdensome requirements.
To this end, we propose a software framework
that transforms a standard fully assembled discretization matrix 
(that might be produced with any standard finite element software) into a series
of structured matrices.
Of course, the underlying mesh used with
the finite element software must coincide with a series of structured
regions (e.g., as in Figure~\ref{PBO:fig:meshRefine}). Additionally, the 
finite element software must provide some minimal information about the 
underlying structured region layout. 

An overall semi-structured solver is being developed within the Trilinos framework\footnote{https://trilinos.github.io} in conjunction with the Trilinos/{\muelu}~\cite{MueLuURL,BergerVergiat2019a} multigrid package.
This solver is not oriented toward matrix-free representations in favor of greater generality,
though some matrix-free performance/memory benefits are sacrificed.
The ideas described in this paper are intended to facilitate the use
of semi-structured solvers within the finite element community and to ultimately
provide significant performance gains over existing fully unstructured 
algebraic multigrid solvers (such as those provided by {\muelu}). 
Section~\ref{sec:StructuredGrids} motivates and describes some semi-structured mesh
scenarios. Section~\ref{sec:RegionVcycle} is the heart of the mathematical
framework, describing the key kernels and their equivalence 
to a standard composite grid multigrid scheme.  Here, the {\vcycle}
application relies heavily on developing a matrix-vector product suitable for matrices
stored in a region-oriented fashion. We also detail 
the hierarchy setup, focusing on the construction of
region-oriented matrices to represent grid transfers and the coarse
discretization matrix.  
Section~\ref{sec:nonInvasive} describes
the framework and the non-invasive application requirements
while 
\Secref{sec:StructuredUnstructured} discusses unstructured regions focusing on the 
treatment of multigrid transfer operators at region interfaces.
We conclude with some numerical experiments to highlight the potential
of such a semi-structured multigrid solver.

\section{Semi-structured grids and mesh abstractions}
\label{sec:StructuredGrids}

Unstructured meshes facilitate the modeling of complex features, but induce 
performance challenges.  Our goal is to provide additional mechanisms to 
address unstructured calculations while furnishing enough structure to reap 
performance benefits.  Our framework centers around block structured meshes
(BSMs).
In our context, it is
motivated by an  existing Sandia hypersonic flow capability where 
the solution quality obtained with block structured meshes 
is noticeably
superior than solutions obtained with fully unstructured meshes\footnote{
This is due to the discretization characteristics 
and mesh alignment with the flying object and with the bow shock.}.
\begin{figure}[ht!]
  \centering
\includegraphics[trim = 1.0in 7.5in 1.0in 6.0in, clip, height = 4.cm,width = 8.cm]{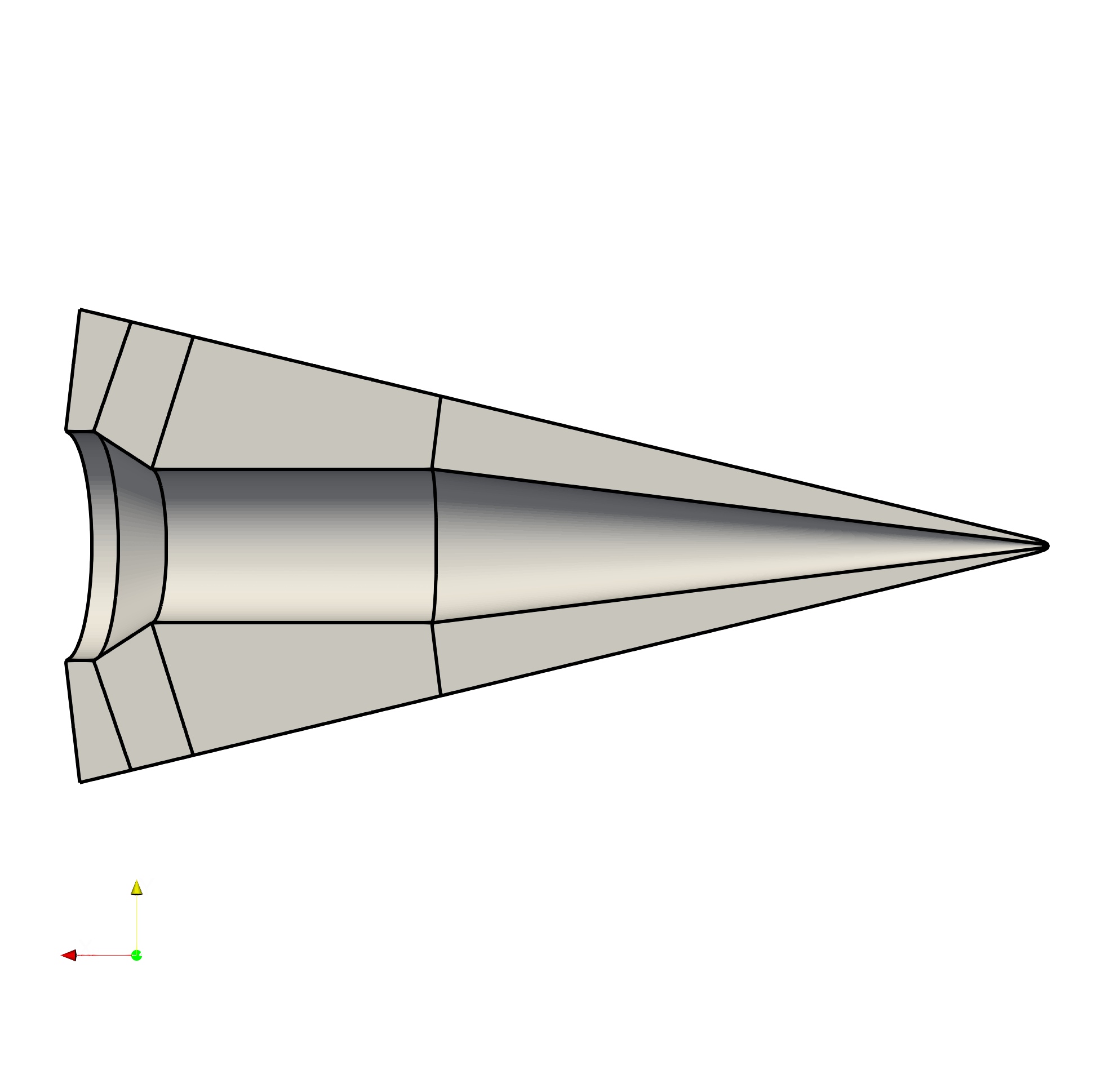}
~~~~~~~
\includegraphics[trim = 5.0in 3.5in 2.3in 1.5in, clip, height = 4.cm,width = 4.cm]{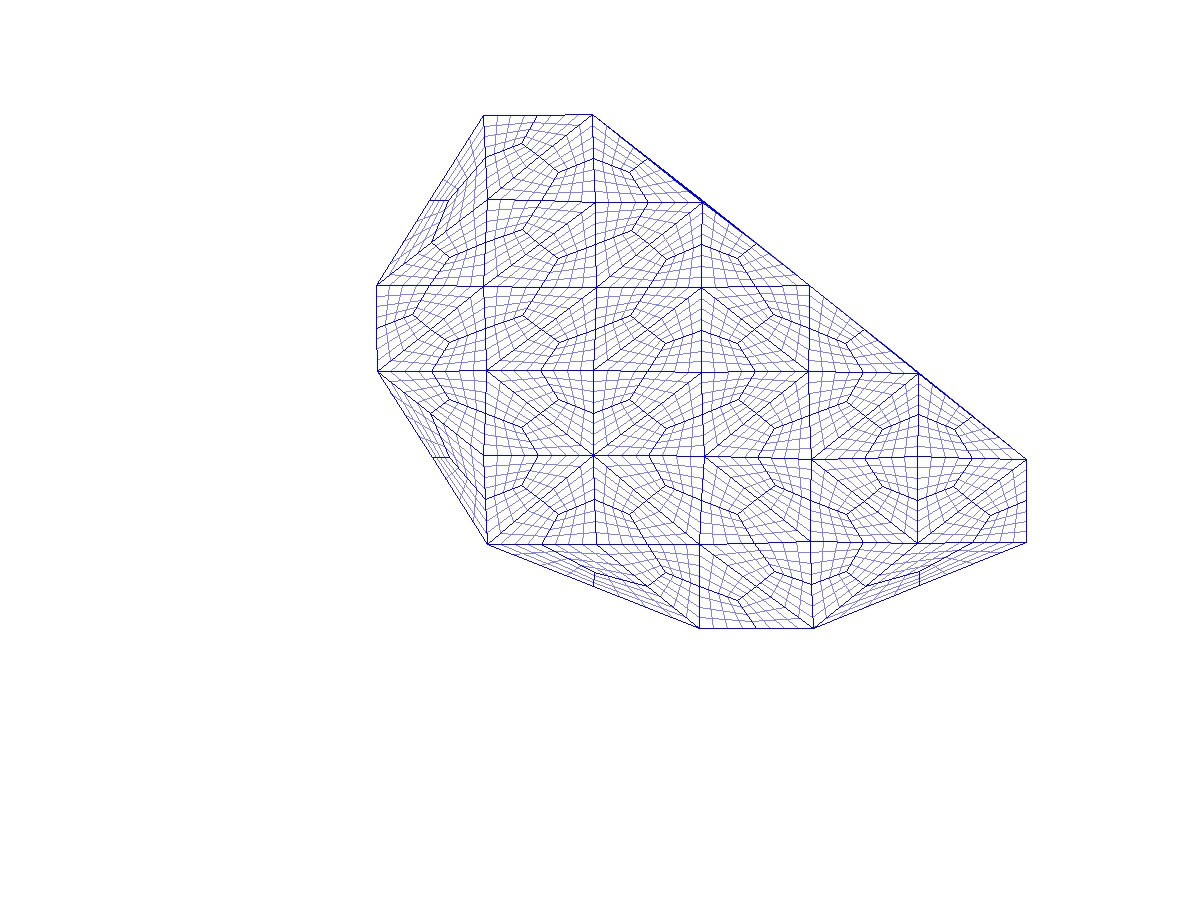}
 \caption{Hypersonic BSM domain (outline of region boundaries depicted; structured grid lines not shown) and BSM/HHG mesh. \label{HHG sample}}
\end{figure}
In this case, BSMs generated by meshing separate components are of significantly greater interest
than meshes of the HHG variety.
Figure~\ref{HHG sample} illustrates a general BSM and a BSM/HHG mesh.

\begin{figure}[ht!]
\centering
\hfill
\includegraphics[trim = 0.5in 0.0in 0.0in 0.0in, clip, height = 3.5cm,width = 3cm]{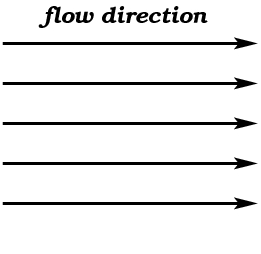}
\includegraphics[trim = 3.1in 0.0in 0.3in 0.1in, clip, height = 3.9cm,width = 6cm]{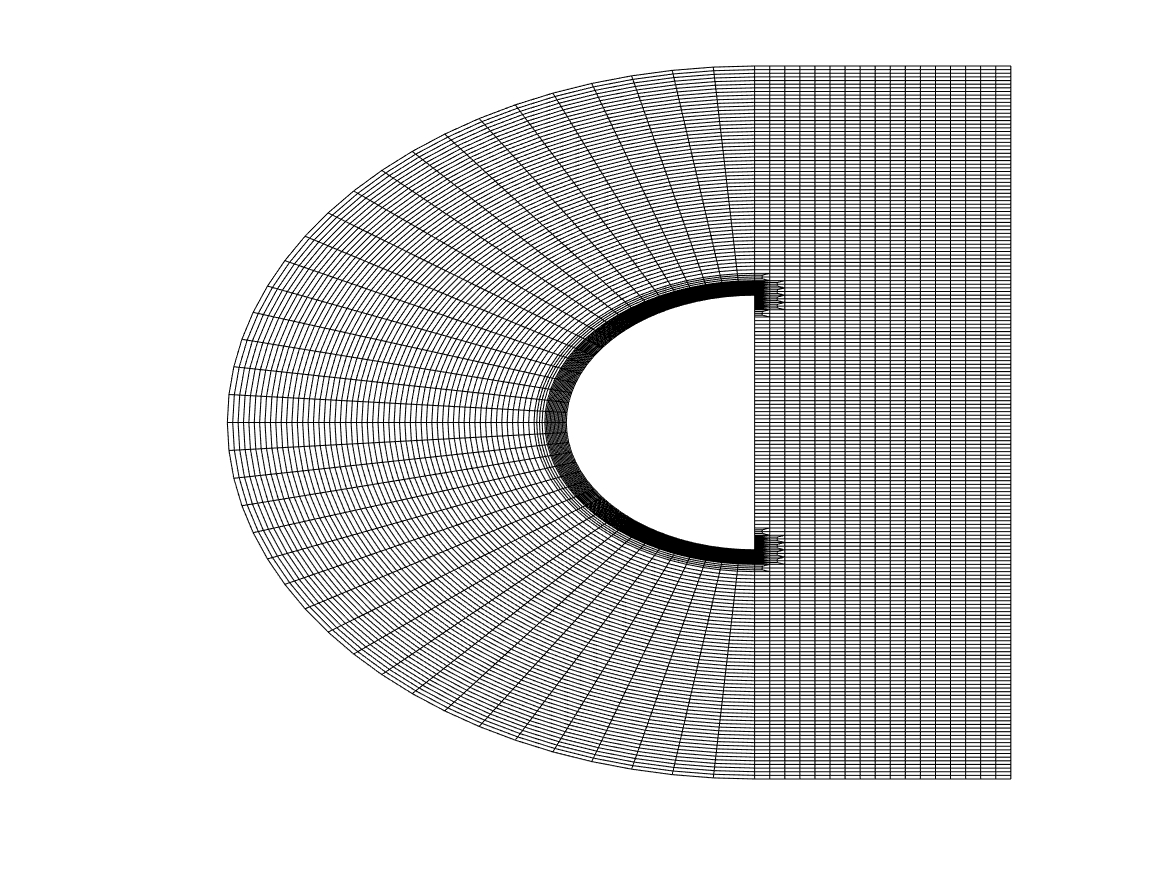}
\hfill
\includegraphics[trim = 1.2in 0.0in 0.3in 0.1in, clip, height = 3.9cm,width = 3cm]{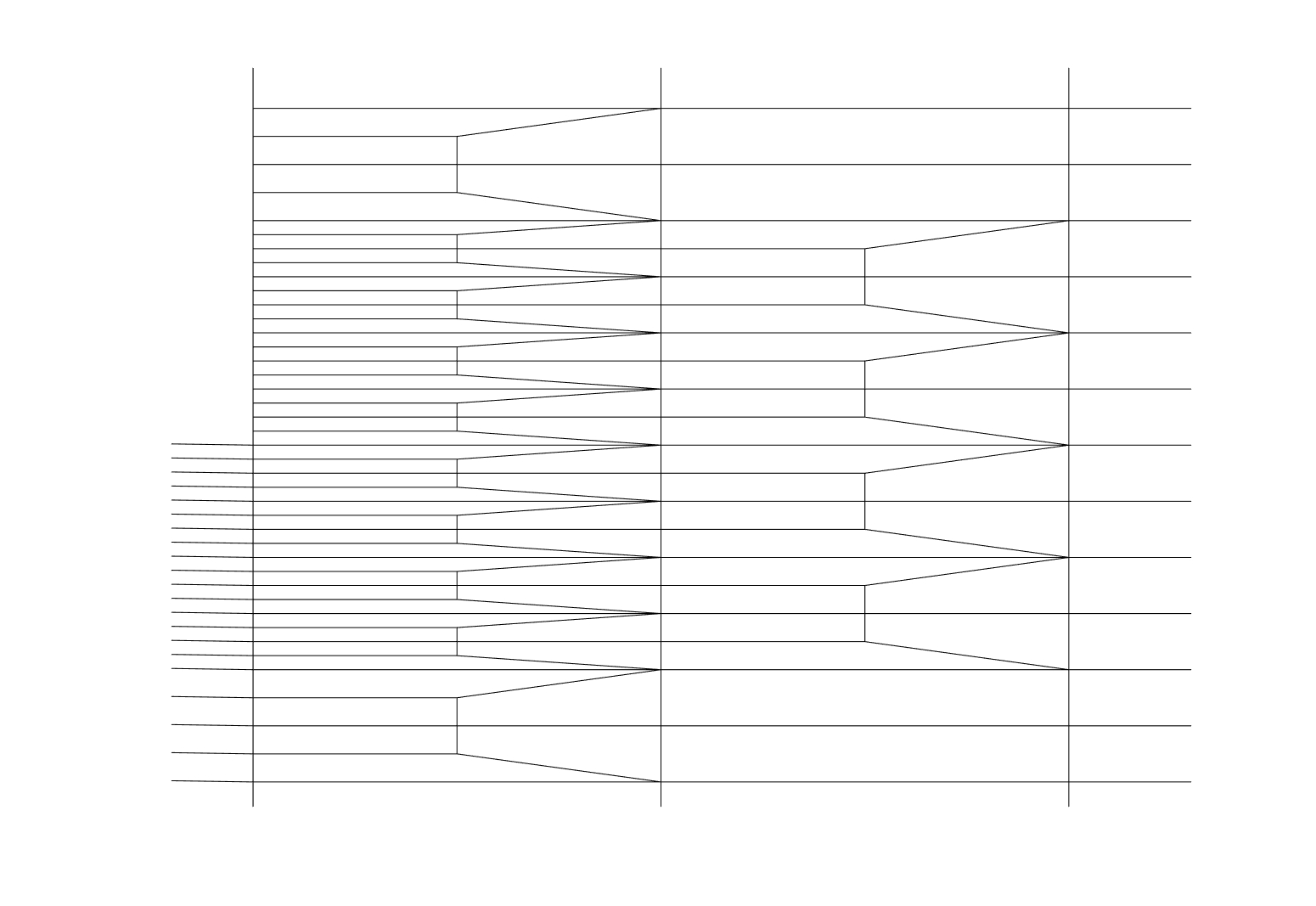}
\hfill~
\caption{Primarily structured mesh with small unstructured regions (left) with a close up view of one of the unstructured regions (right). \label{cyl}}
\end{figure}

While BSMs provide a certain degree of flexibility, unstructured meshes are
often natural to capture complex features locally.
Figure~\ref{radtri} illustrates some scenarios where unstructured regions
might be desirable. Figure~\ref{cyl} shows another case which is
similar to our motivating/target hypersonic example. 
In our hypersonic
problem, refined structured meshes are needed in sub-domains upstream of the
obstacle. In the wake area, however, much lower resolution meshes
(and unstructured meshes) can be employed. In this case, unstructured 
mesh regions can be used to  
transition between structured meshes where modeling characteristics allow for
a large difference in resolutions.
Specifically, two conformal structured meshes could have been
used to represent the domain in Figure~\ref{cyl} (one upstream and the other in the wake). However, the use 
of small unstructured mesh regions allows for a much coarser version of the wake mesh, even though most of the
wake can still be represented with structured mesh regions. 

Our ultimate target is a mesh that includes an arbitrary number of 
structured or unstructured regions that conform at region interfaces.
In this ideal setting, a finite element practitioner would have 
complete freedom to decide the layout of the mesh regions that is 
most suitable for the application of interest. Of course, such a mesh
must be suitably partitioned over processors so that the structured
regions can take advantage of structured algorithms and that the 
overall calculation is load balanced. Here, load balance must take
into account that calculations in unstructured regions will likely
be less efficient than those in structured regions.
While our framework has been designed with this ultimate target in
mind, some aspects of the present implementation limit the
current software to the restriction of one region per processor.


\section{Region-oriented multigrid} \label{sec:RegionVcycle}
We sketch the main ideas behind a region-oriented version of a multigrid solver. 
In some cases, this region-oriented multigrid is mathematically identical to a classical multigrid solver, 
though implementation of the underlying kernels will be different. 
In other cases, it is natural to introduce modest numerical changes to the region-oriented version 
(e.g., a region-local {\GS} smoother).  To simplify notation, we describe only a two level
multigrid algorithm, as the extension to the multilevel case is straight-forward.
\REMOVE {
Algorithm~\ref{alg:RegionVCycle} illustrates a standard multigrid {\vcycle} with our notation.
\begin{algorithm}
\caption{Multigrid {\vcycle} to solve $\MD{A}_\ell \VD{u} = \VD{b}$ .\label{alg:RegionVCycle}}
\begin{algorithmic}
\STATE $\text{function}\;\textsf{MGV}(\MD{A}_{\ell}, \VD{u}, \VD{b},\ell):$
\IF {$ \ell \ne \ell_\text{max}  $}
       \STATE $\VD{u} \gets  \mathcal{S}_\ell (\MD{A}_\ell , \VD{u}, \VD{b}) $
       \STATE $\VD{r} \gets \VD{b} - \MD{A}_{\ell} \VD{u} $
       \STATE $\VD{u}_c \gets 0 $
       \STATE $\VD{u}_c \gets \textsf{MGV}(\MD{A}_{\ell+1},\VD{u}_c, \MD{R}_\ell \VD{r}, \ell\hskip-3pt+\hskip-3pt1)$
       \STATE $\VD{u} \gets \VD{u} + \MD{P}_\ell \VD{u}_c $
\ELSE
       \STATE $\VD{u} \gets \MD{A}_{\ell}^{-1} \VD{b} $
\ENDIF
\end{algorithmic}
\end{algorithm}
}
Figure~\ref{standard MG cycle} provides a high-level illustration of the setup and solve phases
of a classical two level multigrid algorithm. 
\begin{figure}
\vskip .1in
\begin{tabbing}
\hskip .4in \= \hskip .16in \= \hskip 2.3in \= ~~~~~\= \kill
\>$\text{function}\;\textsf{mgSetup}\MD(A,\Psi)$ \>\> $\text{function}\;\textsf{mgCycle}(\MD{A}, \VD{u}, \VD{b}):$ \\[3pt]
\>\> $\textsf{sData} \gets \textsf{smootherSetup}\MD(A)$ \>\>     $\VD{u} \gets  \mathcal{S} (\MD{A}, \VD{u}, \VD{b},\textsf{sData}) $ \\
\>\> $P~~~~~ \gets \textsf{construct}P(\MD{A})$ \>\> $\VD{r} \gets \VD{b} - \MD{A} \VD{u} $ \\
\>\> $R~~~~~ \gets  P^T $ \>\> $\VD{\bar{u}} \gets 0 $ \\
\>\> $\MD{\bar{A}}~~~~~ \gets R A P $ \>\> $\VD{\bar{u}} \gets \textsf{solve}(\MD{\bar{A} },\VD{\bar{u}}, \MD{R} \VD{r})$ \\
\>\>\>\>$\VD{u} \gets \VD{u} + \MD{P} \VD{\bar{u}}$ \\
\end{tabbing}
\caption{Two level multigrid for the solution of $\MD{A} \VD{u} = \VD{b}$.\label{standard MG cycle}}
\end{figure}
%
%
\REMOVE {
\begin{algorithm}
\caption{Region {\multigrid} {\vcycle}}
\label{alg:RegionVCycle}
\DontPrintSemicolon
\SetAlgoLined
\emph{{\Multigrid} algorithm, this is a bit messed up}\;
\Begin{
\emph{Check~$\indLevel$ for coarsest level}\;
\If{$\indLevel == \numLevels - 1$}{
	Form composite coarse level operator~$\levelL{\HHGComposite{\linMat}}$\;
	Solve for the coarse level problem: $\levelL{\linSol} \leftarrow \inv{\levelL{\linMat}}\levelL{\linRhs}$\;
  \Return $\levelL{\linSol}$\;
}
\emph{Apply {\multigrid} method recursively}\;
\Else{
  \emph{Apply~$\numPreSmooth$ pre-smoothing sweeps using~$\smoother$ as a smoother}\;
  \For{$i\leftarrow 1$ \KwTo $\numPreSmooth$}{
    $\levelL{\linSol}\leftarrow\levelL{\smoother}\left(\levelL{\linSol},\levelL{\linRhs}\right)$\;
  }
  \BlankLine
  \emph{Calculate coarse level residual}\;
  
  \BlankLine
  \emph{Initialize coarse level error}\;
  $\levelL{\error} \leftarrow 0$\;
}
}
\end{algorithm}
}
Therein, $A$ refers to the discretization operator on the fine level of the multigrid hierarchy. ${\cal S}$ denotes the 
fine level multigrid smoother.  $P$ interpolates solutions from the coarse level to the fine level while
$R$ restricts residuals from the fine level to the coarse level. $\textsf{sData}$ refers to any pre-computed 
quantities that might be used in the smoother (e.g., ILU factors). 
Coarse level matrices and vectors are delineated
by over bars (e.g., $\bar{A}$ is the coarse level discretization matrix and $\bar{u}$ is the 
coarse level correction).  In this paper, $R$ is always taken as the transpose of $P$, though
the ideas easily generalize to other choices for $R$. Finally, the 
coarse discretization is defined by the projection 
$$
\bar{A} = R A P .
$$
For a two-level method, $\textsf{solve}()$ might correspond to a direct factorization solution
method or possibly coarse level smoother sweeps. In these cases, $\textsf{mgSetup}()$ must
include the setup of the $LU$ factors or coarse level smoothing data. A multilevel algorithm
is realized by instead defining $\textsf{solve}()$ to be a recursive invocation of $\textsf{mgCycle}()$.

The region-oriented multigrid cycle is identical to this standard cycle. The only differences are that
\begin{itemize}
\item $A, \bar{A}, R,$ and $P$  are stored in a region-oriented format,
\item all vectors (e.g., approximate solutions, residuals) are stored in a region-oriented format,
\item all operations (e.g., smoothing kernels) are implemented in a region-oriented fashion with the
exception of the coarsest direct solve.
\end{itemize}
To describe region-oriented multigrid, we begin with a definition of the region layout for vectors and matrices. 
The creation of region-oriented matrices and vectors is delineated in two parts.
The first part focuses on the hierarchy
construction of region-oriented operators when region-oriented operators are provided on the finest level. The second
part then proposes a mechanism for generating  the finest level region-oriented operators 
using information that a standard finite element application can often supply.

\REMOVE{
We briefly sketch the main ideas associated with a region-oriented version of a multigrid solver. The general goal is that this solver be 
mathematically identical (or almost identical) to standard multigrid, only that the implementation
be different. Specifically, 
}
\subsection{Region matrices and vectors}
\begin{figure}
  \centering
\includegraphics[trim = 2.0in 2.3in 1.6in 2.3in, clip, height = 3.cm,width = 9.cm]{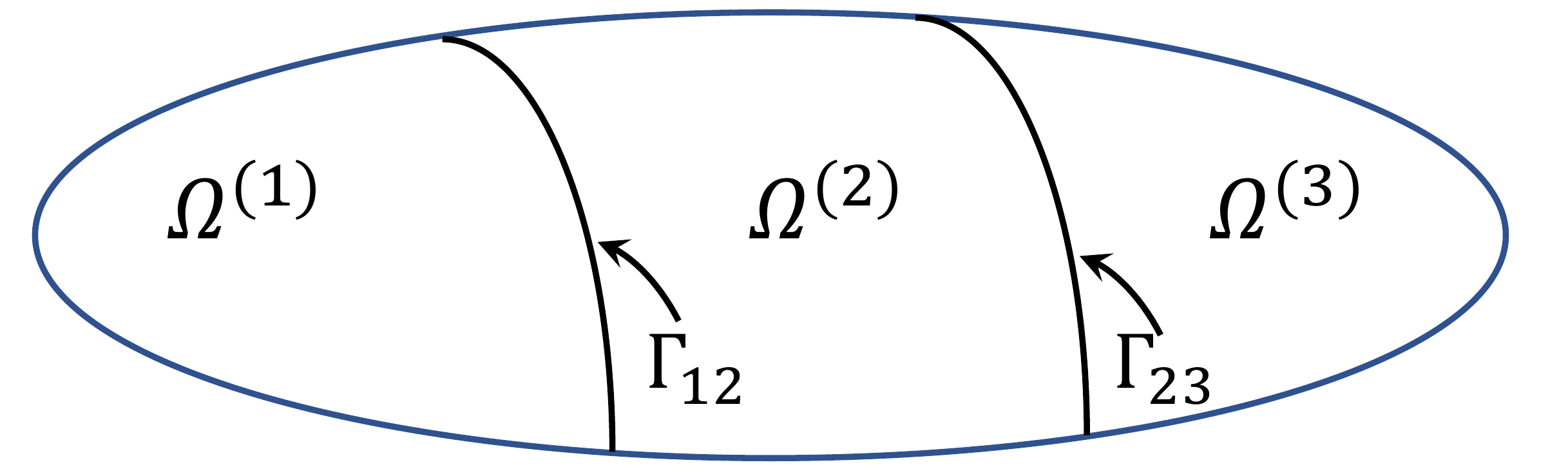}
\caption{Sample domain decomposed into three sub-regions.\label{domains}}
\end{figure}

Consider the discretization of a partial differential equation (PDE) and boundary conditions on a 
domain $\dom$ resulting in the discrete matrix problem
$$
\HHGComposite{\linMat} u = b  .
$$
Often we will refer to the $n \times n$ matrix $\HHGComposite{\linMat} $ as the composite
matrix. Consider now a decomposition of the domain $\dom$ into a set of $\nRegions$ sub-regions 
$\region{\dom}{i}$ such that 
$$
\dom = \cup_{i=1}^{\nRegions}~\region{\dom}{i} .
$$
These regions only overlap at interfaces where they meet (e.g., see Figure~\ref{domains}).
That is, 
$$
\Boundss{i}{j} = \Boundss{j}{i} = \region{\dom}{i}\cap\region{\dom}{j} .
$$
In general, several regions might also meet at so-called corner vertices.
The regions can now be used to split the composite matrix 
such that
\begin{equation} \label{extraction sum}
\HHGComposite{\linMat} = \sum_{1 \le k \le \nRegions} \region{\linMat}{k}
\end{equation}
where 
\begin{equation} \label{extraction restriction1}
\region{\linMat}{k}_{ij}
\ne 0 \qquad \Rightarrow \qquad i,j \in S^{(k)}.
\end{equation}
and 
\begin{equation} \label{extraction restriction2}
\region{\linMat}{k}_{ij}
\ne 0 \qquad \Rightarrow \qquad \HHGComposite{\linMat_{ij}} \ne 0 .
\end{equation}
Here, $S^{(k)}$ is the set of mesh nodes located within $\region{\dom}{k}$ (including those on the interface).
While formally $\region{\linMat}{k}$ is $n \times n$, most rows are identically zero (i.e., rows not
associated with $ S_k$) and so the associated software would only store or compute on non-zero rows.

Mathematically, a region vector is an extended version of a composite vector that 
we express as
$$
\HHGRepRegional{v}^T = \left[~\HHGRepRegional{v}_1^T , ~~~...,~~~ \HHGRepRegional{v}_\nRegions^T ~\right]^T 
$$
where double brackets denote regional representations,
$v$ is the associated composite vector, and $\HHGRepRegional{v}_k$ is a sub-vector of $\HHGRepRegional{v}$ that
consists of all degrees-of-freedom (dofs) that are co-located with the composite dofs given by $S^{(k)}$.
We assume without loss of generality that region dofs within the same region are ordered consecutively
(because region dofs can be ordered arbitrarily).
As composite interface dofs
reside within several regions, the vector $\HHGRepRegional{v}$ will be of length $n_r$ where $n_r \geq n$. 
If we consider a scalar problem and discrete representation of the example given in Figure~\ref{domains},
$\HHGRepRegional{v}$ consists of two dofs for each composite dof on
$\Boundss{1}{2}$ and $\Boundss{2}{3}$.

A region framework can now be understood via a set of boolean transformation matrices.
In particular, a composite vector must be transformed to a region vector where
dofs associated with interfaces are replicated.
To do this, consider an $n \times n_r$ boolean matrix
that maps regional dofs to composite dofs.
Specifically, a nonzero in the $i^{th}$ row and $j^{th}$ column implies that the  $j^{th}$ regional unknown is
co-located with the $i^{th}$ composite unknown.  Each column of $\Psi$ has only
one non-zero entry while the number of non-zeros in a row $i$ of $\Psi$ is equal
to the number of regions that share the $i^{th}$ composite dof.
Thus, a composite
vector $\HHGComposite{v}$ is mapped to a region vector $\HHGRepRegional{v}$ via
$\HHGRepRegional{v} = \Psi^T  \HHGComposite{v}.$
The following properties are easily verified:
\vskip .1in
\begin{tabbing}
\hskip .06in \= \hskip 1.4in \= ~~~~~\= \kill
\> $\Psi \Psi^T$                     \> is a diagonal matrix where the $(j,j)$ entry is the number of region dofs that are\\ 
                                   \>\> co-located with the $j^{th}$ composite dof; \\[3pt]
\>$w = \Psi  \HHGRepRegional{v}$                      \> defines the $j^{th}$ element of $w$ as the sum of the co-located regional elements in\\
                                   \>\> $v$ associated  with composite dof $j$; \\[3pt]
\>$\HHGRepRegional{w}=\Psi^T \Psi \HHGRepRegional{v}$                  \> defines the $j^{th}$ element of $w$ as the sum of the co-located regional elements in\\
               \>\> $v$ associated with regional dof $j$; \\[3pt]
\>$w=(\Psi \Psi^T)^{-1} \Psi \HHGRepRegional{v}$      \>  defines the $j^{th}$ element of $w$ as the average of the co-located regional elements in \\
                                   \>\> $v$ associated with composite dof $j$; \\[3pt]
\>$\HHGRepRegional{w}\hskip -.04in =\hskip -.04in \Psi^T (\Psi \Psi^T)^{-1}\Psi \HHGRepRegional{v}$\> defines the $j^{th}$ element of $w$ as the average of the co-located regional elements in\\
                                   \>\> $v$ associated with regional dof $j$.
\end{tabbing}
Further, one can partition the columns of $\Psi$ in a region-wise fashion such that
\begin{equation} \label{psi partition}
\Psi = \left[\Psi_1 , ~~~...,~~~ \Psi_\nRegions \right].
\end{equation}
Thus, $\Psi_k^T$ maps composite dofs to only region $k$'s dofs,
i.e., $\HHGRepRegional{v}_k = \Psi_k^T v $.
The following additional properties hold:
\vskip .1in
\begin{tabbing}
\hskip .2in \= \hskip 1.25in \= ~~~~~\= \kill
\> $\Psi_k \Psi_k^{T}$   \> filters out dofs not associated with region $k$. In particular, $\Psi_k \Psi_k^{T}$ maps region \\
                           \>\> vectors to new region vectors where the only nonzero matrix entries correspond \\
                           \>\> to an identity block for dofs associated with region $k$; \\[3pt]
\> $S = \Psi_k \Psi_k^T S$ \> if and only if $S$ only contains nonzeros in rows associated with region $k$;  \\[3pt]
\> $S = S \Psi_k \Psi_k^T$ \> if and only if $S$ only contains nonzeros in columns associated with region $k$;  \\[3pt]
\> $\Psi_k^{T} S \Psi_k$   \> is the submatrix of $S$ corresponding to the rows and columns of region $k$.
\end{tabbing}
\vskip .15in
\noindent
The boolean transformation matrices are not explicitly stored/manipulated in our software. Instead,
functions are implemented to perform some of the properties listed above (e.g., averaging
interface values).

A block diagonal region matrix can now be defined as
\begin{equation} \label{block diagonal region matrix}
\HHGRegional{A} =
\begin{pmatrix}
\Psi_1^T A^{(1)} \Psi_1& & & & \cr
                               &.& & & \cr
                               & &.& & \cr
                               & & &.& \cr 
                               & & & & \Psi_\nRegions^T A^{(\nRegions)} \Psi_\nRegions
\end{pmatrix}
.
\end{equation}
Here, we employ a slightly different bracket symbol to emphasize that rows/columns associated with
co-located dofs do not necessarily have the same values in this regional representation.

\begin{lemma} \label{Ar Acomp relationship}
Let $\HHGRegional{A}$ be defined by \eqref{block diagonal region matrix} and $\Psi$ be the boolean 
transformation matrix between region dofs and vector dofs. Then,
\begin{equation} \label{A equivalence}
\Psi \HHGRegional{A} \Psi^T  = A 
\end{equation}
when each split matrix $A^{(k)}$ only contains nonzeros in rows and columns associated
with region $k$'s dofs.
\end{lemma}

\begin{proof}
\begin{align}
\Psi \HHGRegional{A} \Psi^T  &= \Psi_1 \Psi_1^T A^{(1)} \Psi_1 \Psi_1^T + ... + 
                       \Psi_\nRegions \Psi_\nRegions^T A^{(\nRegions)} \Psi_\nRegions \Psi_\nRegions^T \label{orig} \\
                &=                 A^{(1)} \Psi_1 \Psi_1^T + ... + 
                                               A^{(\nRegions)} \Psi_\nRegions \Psi_\nRegions^T  \label{simp1}\\
                &=                 A^{(1)}                 + ... + 
                                               A^{(\nRegions)}       \label{simp2}\\
                &=                 A
\end{align}
where the simplifications to obtain \eqref{simp1} and \eqref{simp2} require that $A^{(k)}$ only have nonzeros in  rows and columns associated with region $k$.
\end{proof}

To rewrite a multigrid {\vcycle} in a region oriented fashion, operations such as matrix-vector
products must be performed with region matrices. For example, matrix-vector products
with the discretization operator in the original multigrid cycle can instead be accomplished
using \eqref{A equivalence}. We also need to replace matrix-vector products associated
with the grid transfers.  For grid transfers, we prefer a different type of region matrix that 
we refer to as replicated interface matrices. Specifically, the replicated interface matrix for 
interpolation is defined by 
\begin{equation} \label{region interp}
\HHGRepRegional{P} = 
\begin{pmatrix}
\Psi_1^T P \cPsi_1& & & & \cr
                               &.& & & \cr
                               & &.& & \cr
                               & & &.& \cr 
                               & & & & \Psi_\nRegions^T P \cPsi_\nRegions 
\end{pmatrix}
\end{equation}
where $\cPsi$ is the boolean matrix associated with the regional to composite transformation on the coarse grid.
Contrary to the standard region matrices, the composite operator (instead of split matrices)
is injected to each of the regions. 
This implies that along the inter-region
interfaces, matrix entries are replicated. 
\begin{lemma} \label{P commute lemma}
\begin{equation} \label{Pcommute}
\HHGRepRegional{P} \cPsi^T = \Psi^T P
\end{equation}
when rows in the matrix $P$  do not contain nonzeros associated with multiple region
interiors (i.e., non-interface dofs from multiple regions).
\end{lemma}
\begin{proof}
\begin{align} 
\HHGRepRegional{P}  \cPsi^T &= \begin{pmatrix} \Psi_1^T P \cPsi_1 \cPsi_1^T \cr
           . \cr
           . \cr 
           . \cr
           \Psi_\nRegions^T P \cPsi_\nRegions \cPsi_\nRegions^T 
\end{pmatrix}
 = \begin{pmatrix} \Psi_1^T P \cr
           . \cr
           . \cr 
           . \cr
           \Psi_\nRegions^T P 
\end{pmatrix}
= \Psi^T P 
\end{align} 
where we use the fact that the matrix $\Psi_k P$ only contains rows associated with 
region $k$ and that this submatrix contains only nonzeros in columns associated with region $k$
(under the assumption that $P$'s rows do not cross multiple region interiors).
\end{proof}
\begin{lemma} \label{R commute lemma}
\begin{equation}  \label{Rcommute}
\cPsi \HHGRepRegional{R} = R \Psi
\end{equation}
when 
\begin{equation} \label{region restriction}
\HHGRepRegional{R} = 
\begin{pmatrix}
\cPsi_1^T R \Psi_1& & & & \cr
                               &.& & & \cr
                               & &.& & \cr
                               & & &.& \cr 
                               & & & & \cPsi_\nRegions^T R \Psi_\nRegions 
\end{pmatrix}
\end{equation}
and $R$ contains no columns where the nonzeros are associated with multiple region interiors.
\end{lemma}
\begin{proof}
Proof omitted as it is essentially identical to the proof for Lemma~\ref{P commute lemma}.
\end{proof}
\vskip .1in
\begin{theorem} \label{rap theorem}
$\cPsi \HHGRepRegional{R} \HHGRegional{A} \HHGRepRegional{P} \cPsi^T  = R A P $
\end{theorem}
\begin{proof}
Follows as a direct result of applying \eqref{Rcommute},  \eqref{Pcommute}, and \eqref{A equivalence}.
\end{proof}

\vskip .1in
Having established basic relationships between region and composite
operations, we now re-formulate the multigrid algorithm primarily in
terms of regional matrices and vectors. This re-formulation must be applied to both the 
multigrid setup phase and the multigrid cycle phase.

\subsection{Multigrid Setup} \label{mg setup}
The multigrid method requires that the discretization matrices, smoothers, and grid transfers
be defined for all levels. For now, let us assume that we have $\Psi$ and $\HHGRegional{A} $
on the finest level. For a two level multigrid method, we must define
$\HHGRepRegional{P},  \HHGRepRegional{R}, \cPsi,$  the regional coarse discretization operator~$\HHGRepRegional{\bar{A}}$,
and the region-based smoothers. 
For grid transfers, we directly create regional forms and never directly form
the composite representation. That is, the composite $P$ and $R$ are only defined
implicitly. In constructing region grid transfers, it is desirable to
leverage standard structured mesh multigrid software\footnote{By ``structured multigrid'', we refer to projection-based multigrid to form coarse operators, but simultaneously exploiting grid structure in the (fine level) discretization. This contrasts geometric multigrid, where coarse levels are formed by an actual re-discretization of the operator on a coarser mesh.} (e.g., apply
structured multigrid software to each region without knowledge of other regions).
However, when creating the regional grid transfers,  the implicitly defined
composite interpolation must not contain any row where different nonzeros are associated with 
different region interiors. Further, stencils from different region blocks
(of the block diagonal interpolation matrix) must be identical for co-located dofs. 
These requirements imply that fine interface vertices must interpolate only from
coarse interface vertices and that interpolation coefficients for fine interface 
dofs have to be identical from neighboring regions.  To satisfy these requirements, we
use standard software in conjunction with some 
post-processing.
In particular, the \emph{standard} grid transfer algorithm must generate some coarse 
points on its region boundary (i.e., the interface) that can be used to fully 
interpolate to
fine vertices on its region boundary.  This is relatively natural for structured
mesh multigrid software.  It is also natural that interpolation stencils 
match along interfaces when using structured multigrid based on linear interpolation
within neighboring regions.  
In this case, grid transfers can be constructed without any communication assuming
that each processor owns one region. That is, each processor constructs the
identical interpolation operator along the interface assuming that each 
processor has a copy of the coordinates and employs the same coarse grid points.
However, if an algorithm is employed
that does not produce identical interpolation coefficients from different
regions, then a natural possibility would be to average the different
interpolation stencils on a shared interface to redefine matching interpolation
stencils at all co-located vertices. This averaging would incur some 
communication when each region is assigned to a different processor.
This type of averaging might be employed if, 
for example, black box multigrid~\cite{Dendy2010a} is used to generate
interpolation within each region as opposed to structured multigrid.
In this way, the region interpolation algorithm will implicitly define a
composite grid interpolation matrix that satisfies \eqref{region interp}.
Regional restriction matrices are obtained by taking the transpose of the
regional interpolation matrices.

Coarse level discretizations can be constructed trivially.  As indicated by Theorem~\ref{rap theorem},
the regional coarse discretization is given by 
\begin{align}
\label{eq:RegionalGalerkinProduct}
\HHGRepRegional{\bar{A}} = \HHGRepRegional{R} \HHGRegional{A} \HHGRepRegional{P} ,
\end{align}
which corresponds to performing a separate triple-matrix product for each diagonal block associated
with each region. When a single region is owned by a single processor, 
no communication is needed in projecting the fine level regional discretization
operator to the coarser
levels.  Given the major scaling challenges of these matrix-matrix operations
within standard AMG algorithms, the importance of being able to
perform this operation in a completely region-local fashion is
significant.  It should be noted, however, that a composite discretization matrix
might be needed at the coarsest level for third-party 
software packages used to provide direct solvers or to further coarsen meshes in an unstructured AMG fashion.
Of course, these composite matrices will only be needed at fairly coarse resolutions
and they can be formed on the targeted level only (i.e., they do not have to be carried through all 
hierarchy levels). Thus, the costs associated with this construction 
via \eqref{A equivalence}  should be modest. 

To complete the multigrid setup, smoothers may require some setup phase.
For Jacobi, {\GS}, and {\Cheby} smoothing,
the diagonal of the composite matrix must be computed during the setup phase. This is easily
accomplished by storing the diagonal of the regional discretization matrix as a regional vector,
e.g.  $\HHGRepRegional{v} = diag(\HHGRepRegional{A})$ using Matlab notation, 
and then simply applying the transformation, i.e., $\Psi^T \Psi \HHGRepRegional{v}$.
For more sophisticated smoothers, it is natural to generate region analogs 
that are not completely equivalent to the composite versions. For example, one 
can generate region-local versions of {\GS} smoothers and Schwarz type methods
where again $\Psi^T \Psi$ may be used to perform sums of nonzeros from different 
regions associated with co-located vertices.
In this paper, we consider Jacobi, {\GS}, and {\Cheby} smoothers.
Some discussion of more sophisticated smoothers can be found in~\cite{Bergen2004a}.

Finally, construction of a coarse level composite operator~$\bar{A}$ is also trivial. 
In particular, $\cPsi$ is just the submatrix of $\Psi$ corresponding to taking
rows associated with coarse composite vertices and columns associated
with the co-located coarse region vertices. 
Thus, it is convenient if the interpolation algorithm 
also provides a list of coarse vertices, though this can be deduced
from the interpolation matrix (i.e., the vertices associated with rows containing only one nonzero).

Having computed the coarse level operator~$\HHGRepRegional{\bar{A}}$ via the recursive application of~\eqref{eq:RegionalGalerkinProduct},
its composite representation is given as
\begin{align}
\label{eq:CoarseCompositeOperator}
\bar{A} = \cPsi\HHGRepRegional{\bar{A}}.
\end{align}
This corresponds to forming sums of matrix rows that correspond to co-located nodes on region 
interfaces. 

\subsection{Multigrid Cycle}
The multigrid cycle consists primarily of residual calculations, 
restriction, interpolation, and smoother applications. 
The composite residual can be calculated 
with region matrices via
\begin{equation} \label{composite residual}
r = b - A u = b - \Psi \HHGRegional{A} \Psi^T u  .
\end{equation}
Normally, however, one seeks to compute the regional form of the residual
using regional representations of $b$ and $u$ via
\begin{equation} \label{regional residual}
\HHGRepRegional{r} = \HHGRepRegional{b} - \Psi^T \Psi \HHGRegional{A} \HHGRepRegional{u}  ,
\end{equation}
which is derived by pre-multiplying \eqref{composite residual} by $\Psi^T$ and recognizing
that $\HHGRepRegional{r} = \Psi^T r  , \HHGRepRegional{b} = \Psi^T b  ,$ and 
$ \HHGRepRegional{u} = \Psi^T u $. Thus, the only difference with a standard residual
calculation is the interface summation given by $\Psi^T \Psi$.
For interpolation, we seek the regional version of interpolation
\begin{align}
\HHGRepRegional{w} &= \Psi^T P v \\
                   &= \HHGRepRegional{P} \cPsi^T v \\
                   &= \HHGRepRegional{P} \HHGRepRegional{v} 
\end{align}
where we used Lemma~\ref{P commute lemma} to simplify the interpolation expression. 
Thus, the interpolation matrix-vector product is identical to a standard 
matrix-vector product, incurring no inter-region communication.

The region version of the restriction matrix-vector product is a bit more complicated. 
We begin by observing that
\begin{align} 
R &= \cPsi \HHGRepRegional{R} \Psi^T (\Psi \Psi^T)^{-1}\label{R relationship} \\
  &= \cPsi \HHGRepRegional{R} \HHGRepRegional{\Psi \Psi^T}^{-1} \Psi^T  \label{simplified R relationship} .
\end{align}
Lemma~\ref{R commute lemma} can be used to verify \eqref{R relationship}. For \eqref{simplified R relationship},
we define an interface version of $\Psi \Psi^T$ analogous to  \eqref{region interp}
and \eqref{region restriction}. Specifically, the $\HHGRepRegional{\Psi \Psi^T}$ matrix is both diagonal 
and block diagonal where the $k^{th}$ block is given by $\Psi_k^T (\Psi \Psi^T) \Psi_k$. By employing a
commuting relationship (whose proof is omitted as it closely resembles that of Lemma~\ref{P commute lemma}),
one arrives at \eqref{simplified R relationship}. Finally, pre-multiplying $w = R v$ by $\cPsi^T$, substituting 
\eqref{simplified R relationship} for $R$, and recognizing that 
$  \HHGRepRegional{w} = \Psi^T w$ and $  \HHGRepRegional{v} = \Psi^T v$, it can be shown that
the desired matrix-vector product relationship is given by 
$$
\HHGRepRegional{w} = \cPsi^T \cPsi \HHGRepRegional{R} \HHGRepRegional{\Psi \Psi^T}^{-1} \HHGRepRegional{v} .
$$
Thus, the restriction matrix-vector product corresponds to region-local scaling, followed by a region-local
matrix-vector product followed by summation of co-located regional quantities. 

\subsection{Region level smoothers}
\label{sec:RegionLevelSmoothers}

Jacobi smoothing is given by
\begin{align*}
\HHGRepRegional{u} \leftarrow \HHGRepRegional{u} + \omega~\HHGRepRegional{\tilde{D}^{-1}} \HHGRepRegional{r}
\end{align*}
with $\HHGRepRegional{r}$ computed via \eqref{regional residual},
$\omega$ is a damping parameter, and $\HHGRepRegional{\tilde{D}}$ is the diagonal of the
composite operator $\HHGComposite{\linMat}$ stored in regional form (as discussed
in Section~\ref{mg setup}). 

Implementation of a classic {\GS} algorithm always requires some care on parallel 
computers, even when using standard composite operators.  Though a high degree of
concurrency is possible with multi-color versions, these are difficult to develop
efficiently and require communication exchanges for each color on message passing
architectures. Instead, it is logical to adapt region {\GS} using domain 
decomposition ideas (as is typically done for composite operators as well). 
The $K$ sweep {\GS} smoother is summarized in \algref{alg:GSSmoother}.
\begin{algorithm}
\caption{{\GS} smoother for region-type problems}
\label{alg:GSSmoother}
\begin{algorithmic}
\REQUIRE $\omega , \HHGRepRegional{A}, \HHGRepRegional{b}, \HHGRepRegional{\tilde{D}}, \HHGRepRegional{u}$
\FOR{$k = 0, \hdots, K-1$}
\STATE $\HHGRepRegional{\delta} = 0$
\STATE compute $\HHGRepRegional{r}$ via~\eqref{regional residual}
\STATE // for each region ...
\FOR{$\ell = 1, \hdots, m$ }
\FOR{$i = 0, \hdots, N^{(\ell)}$}
\STATE $ r^{(\ell)}_i = r^{(\ell)}_i = - \Sigma_{j} A^{(\ell)}_{ij} \delta^{(\ell)}_j  $
\STATE $ \delta^{(\ell)}_i = \omega r^{(\ell)}_i / \tilde{d}^{(\ell)}_{ii}  $
\STATE $ u^{(\ell)}_i = u^{(\ell)}_i  + \delta^{(\ell)}_{i}  $
\ENDFOR
\ENDFOR
\ENDFOR
\end{algorithmic}
\end{algorithm}
Here, the notation $ r^{(\ell)}_i $ refers to the $i^{th}$ component of the $\ell^{th}$ region's vector while $A^{(\ell)}_{ij}$ refers to a particular nonzero in region $\ell$'s matrix. The intermediate quantity $ \delta^{(\ell)}_i $ is used to update
the local solution and the local residual.  Notice that the only communication 
is embedded 
within the residual calculation at the top of the outer loop. This low 
communication version of the algorithm differs from true {\GS}
in that a region's updated residual only takes into account solution changes 
within the region. This means that solution values along a shared interface are 
not guaranteed to coincide during this state of the algorithm. 

{\Cheby} smoothing relies on optimal {\Cheby} polynomials tailored to reduce errors within the eigenvalue interval~$\lambda_i\in[\eigValMin,\eigValMax]$
with~$\eigValMin$ and~$\eigValMax$ denoting the smallest and largest eigenvalue of interest of the operator~$\HHGRepRegional{A}$.
The largest eigenvalue is obtained by a few iterations of the power method.
Following the {\Cheby} implementation in {\ifpackTwo}~\cite{Prokopenko2016b},
we approximate this interval by $[\eigValMin,\eigValMax]\approx[\alpha,\beta]$
with~$\alpha = \tildeEigValMax / \eigRatio$ and $\beta = \boostFactor\tildeEigValMax$
where~$\tildeEigValMax$ is the estimate obtained via the power method, 
~$\eigRatio $
denotes a ratio that is either user supplied or given by the coarsening rate between levels
(defaulting to $\eigRatio=20$)
and~$\boostFactor$ is the so-called ``boost factor'' (often defaulting to~$\boostFactor=1.1$).
The {\Cheby} smoother up to polynomial degree~$K$ is summarized in \algref{alg:ChebyshevSmoother}.
\begin{algorithm}
\caption{{\Cheby} smoother for region-type problems}
\label{alg:ChebyshevSmoother}
\begin{algorithmic}
\REQUIRE $\theta = \frac{\alpha + \beta}{2}, \delta = \frac{2}{\beta - \alpha}, \HHGRepRegional{A}, \HHGRepRegional{\tilde{D}}, \HHGRepRegional{u}, \HHGRepRegional{r}$
\STATE $\rho = \left(\theta\delta\right)^{-1}$
\STATE $\HHGRepRegional{d} = \frac{1}{\theta}\delta\HHGRepRegional{\tilde{D}^{-1}}\HHGRepRegional{r}$
\FOR{$k = 0, \hdots, K$}
\STATE $\HHGRepRegional{u} = \HHGRepRegional{u} + \HHGRepRegional{d}$
\STATE compute $\HHGRepRegional{r}$ via~\eqref{regional residual}
\STATE $\rho_{\mathrm{old}} = \rho$
\STATE $\rho = \left(2\theta\delta - \rho_{\mathrm{old}}\right)^{-1}$
\STATE $\HHGRepRegional{d} = \rho\rho_{\mathrm{old}}\HHGRepRegional{d} + 2\rho\delta\HHGRepRegional{\tilde{D}^{-1}}\HHGRepRegional{r}$
\ENDFOR
\end{algorithmic}
\end{algorithm}

\subsection{Coarse level solver}
\label{sec:CoarseLevelSolver}
The region hierarchy consists of $\numLevelsRegion$ levels~$\indLevel \in \{0,\hdots,\numLevelsRegion-1\}$.
Having computed the coarse composite operator~$\bar{A}$ via \eqref{eq:CoarseCompositeOperator} on level~$\numLevelsRegion-1$,
we construct a coarse level solver for the region MG hierarchy. We explore two options:
\begin{itemize}
\item {\bf Direct solver:} If tractable, a direct solver relying on the factorization~$\bar{A} = \bar{L} \bar{U}$ is constructed.
As usual, its applicability and performance (especially {\wrt} setup time) largely depend on the number of unknowns on the coarse level.
\item {\bf AMG {\vcycle}:} If~$\bar{A}$ is too large to be tackled by a direct solver,
one can construct a standard AMG hierarchy with an additional $\numLevelsAMG$ levels.
The coarse level solve of the region MG cycle is then replaced by a single {\vcycle} using (SA-)AMG~\cite{Vanek1996a}.
This AMG hierarchy requires only the operator~$\bar{A}$ and its nullspace, which can be extracted from the region hierarchy.
The AMG {\vcycle} itself will create as many levels as needed, such that its coarsest level can be addressed using a direct solver.
The number of additional levels for the AMG {\vcycle} is denoted by~$\numLevelsAMG$.
For efficiency, load re-balancing is crucial.
(Note that the total number of levels is now $\numLevels = \numLevelsRegion + \numLevelsAMG - 1$,
where the subtraction by one reflects the change of data layout from region to composite format without coarsening.)
\end{itemize}
The latter option is also of interest for problems, where the regional fine mesh has been constructed through regular refinement of an unstructured mesh.
Here, the region MG scheme can only coarsen until the original unstructured mesh is recovered.
AMG has to be used for further coarsening.
Assuming one MPI rank per region, i.e. one MPI rank per element in the initial unstructured mesh,
the need for re-balancing (or even multiple re-balancing operations throughout the AMG hierarchy) becomes obvious.

\subsection{Regional multigrid summary}
\label{regionalMG summary}
To summarize, the mathematical foundation and exact equivalence with standard
composite grid multigrid requires that
\begin{enumerate}
\item the composite matrix be split according to \eqref{extraction sum} such 
that each piece only includes nonzeros defined on its corresponding region; 
\item each row (column) of the composite interpolation (restriction) matrix 
      cannot include nonzeros associated with multiple region interiors;
\end{enumerate}
Thus, co-located fine interpolation rows consist only of nonzeros associated 
with coarse co-located vertices. Likewise, co-located coarse restriction 
columns only include nonzeros associated with fine co-located vertices. 
Finally, the grid transfer condition implies that regional forms of 
interpolation (restriction) must have matching rows (columns) associated
with co-located dofs. It is important to notice that if the region interfaces are
not curved or jagged and if linear interpolation is used to define the grid transfer
along region interfaces (where fine interface points only interpolate from coarse
points on the same interface), then each region's block of the block interpolation
operator can be defined independently as long as the selection of coarse points
on the interface match. That is, the resulting region interpolation operator will
satisfy the Lemma conditions without the need for any communication. 
If, however, a more algebraic scheme is used to generate the inter-grid transfers,
then some communication might be needed to ensure that the interpolation operators
satisfy the Lemma conditions at the interface. This would be true if a 
black box multigrid~\cite{Dendy2010a} is used to define the grid transfers
or if a more general algebraic multigrid scheme 
such as smoothed aggregation~\cite{Vanek1996a} is used to define grid transfers.
This is discussed further in Section~\ref{sec:StructuredUnstructured}.

Figure~\ref{regional MG cycle} summarizes the regional version of the two level algorithm.
\begin{figure}
\vskip .1in
\begin{tabbing}
\hskip .4in \= \hskip .16in \= \hskip 2.3in \= ~~~~~\= \kill
\>$\text{function}\;\textsf{mgSetup}(\HHGRegional{A})$ \>\> $\text{function}\;\textsf{mgCycle}(\HHGRegional{A}, \HHGRepRegional{u}, \HHGRepRegional{b}):$ \\[3pt]
\>\> $\HHGRepRegional{D} \gets \textsf{diag}(\Psi^T \Psi~ \textsf{diag}(\HHGRegional{A}))$ \>\>     
$\HHGRepRegional{u} \gets \textsf{applySmoother}(\HHGRepRegional{u}, \HHGRepRegional{b}, \HHGRegional{A})$\\[1pt]

\>\> $\HHGRepRegional{P} \gets \textsf{construct}P(\HHGRegional{A})$ \>\> 
$\HHGRepRegional{r} \gets \HHGRepRegional{b} - \Psi^T \Psi \HHGRegional{A} \HHGRepRegional{u}  $\\
\>\> $\HHGRepRegional{R} \gets  \HHGRepRegional{P}^T $ \>\> $\HHGRepRegional{\bar{u}} \gets 0 $ \\
\>\> $\MD{\HHGRegional{\bar{A}}}\gets \HHGRepRegional{R} \HHGRegional{A} \HHGRepRegional{P} $ \>\> $\HHGRepRegional{\bar{u}} \gets \textsf{solve}(\HHGRegional{\bar{A} },\HHGRepRegional{\bar{u}}, 
\cPsi^T \cPsi \HHGRepRegional{R} \HHGRepRegional{\Psi \Psi^T}^{-1} \HHGRepRegional{r})$\\[1pt]
\>\>$\bar{\Psi} \gets \textsf{inject}(\Psi)$
\>\>$\HHGRepRegional{u} \gets \HHGRepRegional{u} + \HHGRepRegional{P} \HHGRepRegional{\bar{u}}$ \\
\end{tabbing}
\caption{Two level regional multigrid for the solution of $\MD{A} \VD{u} = \VD{b}$.\label{regional MG cycle}}
\end{figure}
Besides the $\textsf{inject}()$ operation, 
the only possible difference during setup is a small modification of $\textsf{construct}P()$ that may be necessary
to ensure that interpolation stencils {\it match} at co-located vertices.
In $\textsf{applySmoother}()$, any region level smoother from \secref{sec:RegionLevelSmoothers} is applied.
The main difference in the $\textsf{solve}()$ phase is the scaling
$\HHGRepRegional{\Psi \Psi^T}^{-1}$, the interface summation $\Psi^T \Psi $, and 
possibly the need to convert between regional and composite forms if third party software is employed
at sufficiently coarse levels.

\section{Non-invasive construction of region application operators} \label{sec:nonInvasive}

To this point, we have assumed that $\Psi$ and $ \HHGRegional{A}$ on the finest 
level are available.  However, most finite element 
software is not organized to generate these. 
Our goal is to limit the burden on application developers by instead
employing a fully
assembled discretization or composite matrix on the finest level. 
In this section, we first
describe the application information that we require to generate $\Psi$.
Then, we describe an automatic matrix splitting or dis-assembly process
so that our software can generate $ \HHGRegional{A}$, effectively via
\eqref{block diagonal region matrix}.

\REMOVE {
In this section,
we describe a matrix dis-assembly process that can transform a fully
assembled matrix~$\HHGComposite{\linMat}$ into something suitable for a region-oriented multigrid
solver, namely the region-oriented operator~$\HHGRegional{\linMat}$. 
More generally, a flexible region-oriented multigrid framework must 
provide an ability to map between region vector/matrix representations and 
composite vector/matrix representations.  Overall, composite matrices/vectors 
are preferred by users and by third-party software packages used to provide direct 
solvers or to further coarsen meshes in an unstructured AMG fashion.
This implies that composite 
matrices/vectors provided
by users must be converted to the regional forms needed by the solver
and then solutions must be converted back to composite forms for the
users. Similarly, when invoking a direct solver or a
conventional AMG package to generate additional levels , regional representations 
must be converted to composite ones and then computed quantities must be converted back
to a regional form for further regional processing.  Additionally, 
a region-oriented framework must be able 
to efficiently apply the region interface operator (or a modified interface 
operator that takes averages).
To do this, some information must be supplied by the application. Our goal
is to limit the burden on application developers by employing a fully
assembled discretization or composite matrix on the finest level. 
}

In addition to fairly standard distributed matrix requirements
(e.g., each processor supplies a subset of owned matrix rows 
and a mapping between local and global indices for the owned
rows), applications must provide information to construct $\Psi$
and to facilitate fast kernels. Specifically, applications
furnish a region id and the number of grid points in each dimension
for regions owned by a processor.
As noted, our software is currently limited in that each processor
owns one entire region. However, we will keep the discussion
general. 

\begin{figure}
\hskip -.17in
  \includegraphics[trim = 0.5in 1.0in 0.0in 0.7in, clip, scale=0.5]{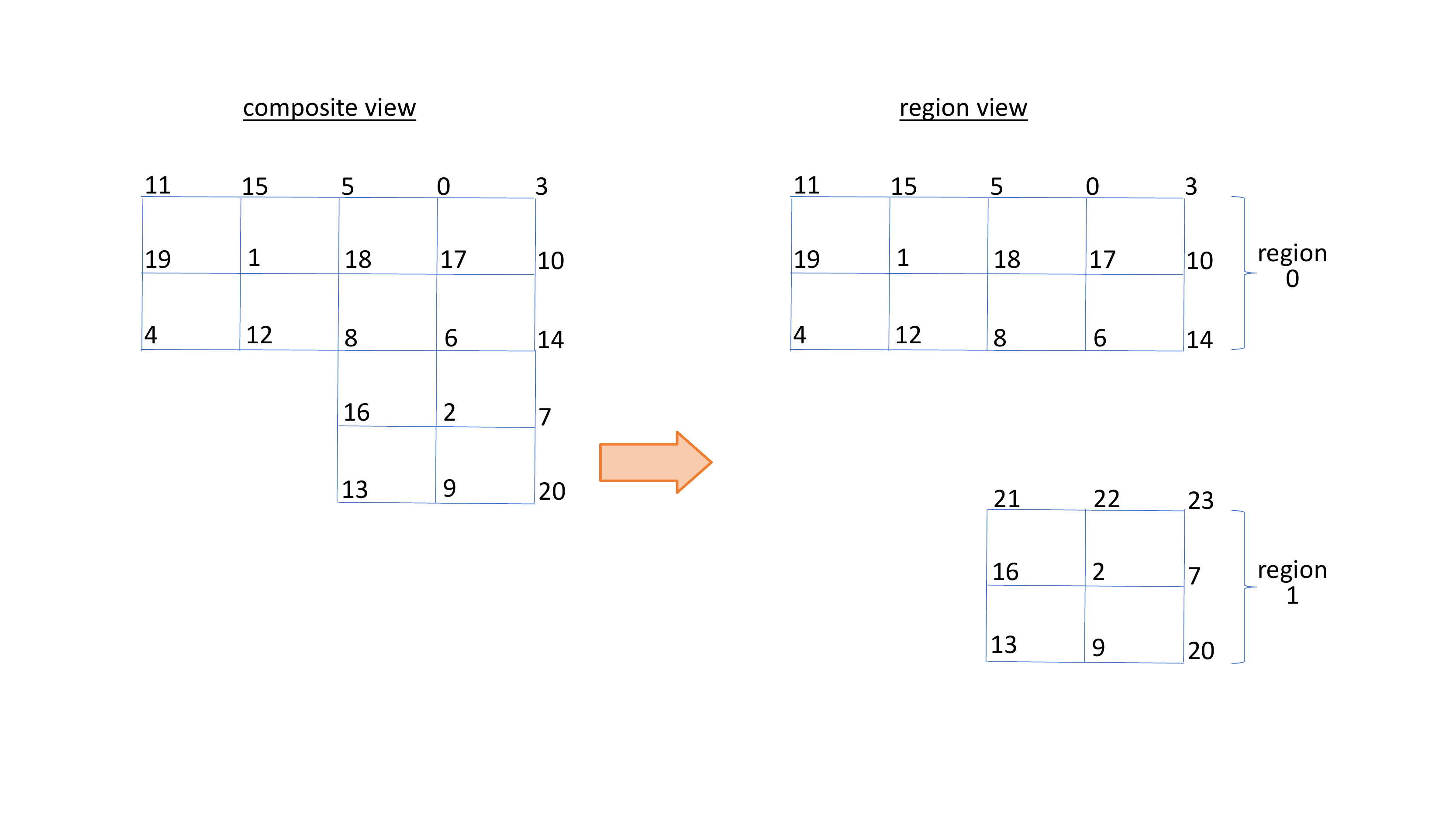}
\caption{Sample user-provided mapping of mesh nodes to regions.\label{regionPerGIDWithGhosts}}
\end{figure}

The main additional requirement is a description of the mesh at the region interfaces.
In particular, it must be known, to which region(s) each node belongs.
If a node is a region-internal node, it only belongs to one region.
If it resides on a region interface, it belongs to multiple regions.
Note that the number of associated regions depends on the spatial dimension,
the location within the mesh, and the region topology.
For example, nodes on inter-region faces (not also on edges and corners),
edges (not also on corners), and corners belong
to 2 regions, 4 regions, and 8 regions respectively
for a three-dimensional problem with cartesian-type cuboid regions.
Figure~\ref{regionPerGIDWithGhosts}
gives a concrete two region example in a two-dimensional setting.
In this example, one processor owns the entire $5 \times 3 $ topmost
rectangular region while another processor owns the bottom most $ 3 \times 2 $
rectangular region. 
The mapping for this example looks as follows:
\begin{itemize}
\item Nodes~$0,1,3,4,5,10,11,12,14,15,17,18,19$ reside in region~$\region{\dom}{0}$.
\item Nodes~$2,7,9,13,16,20$ reside in region~$\region{\dom}{1}$.
\item Nodes~$6,8,14$ are located on the region interface and belong to both regions~$\region{\dom}{0}$ and~$\region{\dom}{1}$.
\end{itemize}
Based on this user-provided mapping data,
we can now ``duplicate'' interface nodes and assign unique GIDs for all replicated interface nodes
and their associated degrees of freedom.
The right-hand side sketch in Figure~\ref{regionPerGIDWithGhosts} illustrates a computed
mapping of global composite ids 
to the global region layout ids. Notice that the only global 
ids to change are the composite ghost ids. Specifically, new global
ids are assigned by the framework to the ghosts associated with the bottom 
processor so that each of the unknowns along a shared interface has a unique
global id.
The overall structured framework can be setup
based on this user-supplied mapping and effectively build the $\Psi$ operator.
Of course, we do not explicitly form $\Psi$,
but build data structures and functions to perform the necessary operations associated with $\Psi$.

To apply \eqref{block diagonal region matrix}, the composite matrix
must first be split so that \eqref{extraction sum}, 
\eqref{extraction restriction1} and \eqref{extraction restriction2} are satisfied.
Mathematically, matrix entries associated with co-located vertices 
must be split or divided between different terms in the summation.
In this paper, we scale any off-diagonal matrix entries by the number of 
regions that share the same edge. 
Formally, scaled entries correspond to
$A_{ij} \ne 0$ such that 
there exist exactly $q ~(\ge 2)~\Psi_k$'s
with a nonzero 
in the $i^{th}$ and $j^{th}$ rows.
If we denote these  $\Psi_k$'s by 
$\Psi_{k_1}, \Psi_{k_2}, ... , \Psi_{k_q}$, then 
$$
A_{ij}^{(k_1)} = 
A_{ij}^{(k_2)} = 
...  = 
A_{ij}^{(k_q)} = 
A_{ij}/q .
$$
The matrix diagonal is then scaled so that the row sum of each region matrix is 
identically zero. 
With $\Psi$ and the splitting choice specified,  the entire multigrid 
cycle is now defined. 
Though this splitting choice is relatively simple, it has no
numerical impact when geometric grid transfers are employed in conjunction
with a Jacobi smoother.  However, some multigrid components such as
region-oriented smoothers (e.g., region-local {\GS}) 
and matrix-dependent algorithms for generating
grid transfers (e.g., black-box multigrid) are affected by the splitting
choice. 
We simply remark
that we have experimented with a variety of scalar PDEs using black-box
multigrid, and this splitting choice generally leads to multigrid convergence
rates that are similar to conventional multigrid algorithms applied
to composite problems. 

While we do not provide the implementation details
associated with computations such as $\Psi_k^T A^{(k)} \Psi_k$ 
and the conversions between regional and composite vectors, 
it is worth pointing out that some implementation aspects 
can leverage ghosting and overlapping Schwarz capabilities found
in many iterative solver frameworks. In our case, some of these 
operations can be performed in a relatively straight-forward fashion 
using Trilinos' import/export mechanism. The import feature is most 
commonly used in Trilinos to perform operations such a matrix-vector 
products.
An import can be used to take vectors without ghost unknowns and create 
a new vector with ghost unknowns obtained from neighboring processors.
This standard import operation is similar to 
transforming a composite vector to 
a region vector. The main difference is that 
only some ghost unknowns (those that correspond to a shared interface)
need to be obtained from neighboring processors.

The import facility is fairly general in that it can 
also be used to replicate matrix rows needed within a standard 
overlapping Schwarz preconditioner. In this case, import takes a 
non-overlapped matrix where each matrix row resides on only one processor
and creates an overlapped matrix, where some matrix
rows are duplicated and reside within more than one sub-domain.
When an overlap of one is used, each
processor receives a  duplicate row for each 
of its ghost unknowns. 
\REMOVE {
These newly received rows will contain
some new column entries that were not previously defined on 
the processor (i.e., not owned by the processor nor part of the
processor's ghost unknowns). These new columns effectively define new ghost 
unknowns that will determine additional rows to be received if 
an overlap of two is desired. Otherwise, the entries in these new columns
will be discarded. Ultimately, a Schwarz procedure seeks to generate
a square subdomain matrix. }
This is similar to the process of generating regional matrices from 
composite matrices
(only requiring rows from a subset of ghosts). 
Once matrix rows (corresponding to interfaces) have been replicated,
they must be modified to satisfy \eqref{extraction sum}. In particular,
any column entries (within interface rows) that correspond to 
connections with neighboring regions must be removed. Further,
entries that have been replicated along the interface must be
scaled in a post-processing step.

In a standard Schwarz preconditioner, solutions obtained on each sub-domain must
be combined. That is, overlapped solution values must be combined (e.g., averaged)
to define a unique non-overlapping solution. For this mapping from overlapped
to non-overlapped, Trilinos contains an export mechanism. This export allows for
different type of operations (e.g., averages or sums) to be used when combining
multiple entries associated with the same non-overlapped unknown.
This is similar to transforming regional vectors to composite vectors.
One somewhat subtle issue is that the
unique region global ids presented in Figure~\ref{regionPerGIDWithGhosts} are not 
needed in an overlapping Schwarz capability, but are needed for the 
region-multigrid framework to perform further operations on the 
region-layout systems.
Thus, the conversions between composite and regional forms has been implemented
in two steps. The first step closely resembles the Schwarz process
and corresponds to the movement of data between overlapped and non-overlapped
representations as just discussed, but without introducing the new global
ids. 
The second step then defines the new global ids
to complete the conversion process.

\section{Structured/unstructured mesh hybrid}
\label{sec:StructuredUnstructured}
We now discuss the adaptation of regional multigrid
to the case where some unstructured regions are introduced
into the grid.  As the mathematical foundation presented earlier
makes no assumptions on grid structure, the requirements summarized 
in Section~\ref{regionalMG summary} still hold. 
The unstructured regions do not introduce
software modifications associated with satisfying the matrix splitting or 
dis-assembly requirements. However, grid transfer construction requires some 
care. In particular,
some pre- and post-processing modifications are needed for the AMG algorithm 
that constructs regional grid transfers within the unstructured regions.
No additional modifications are needed to produce structured grid
multigrid transfers within the structured regions.

\begin{figure}
  \centering
  \includegraphics[scale=0.3]{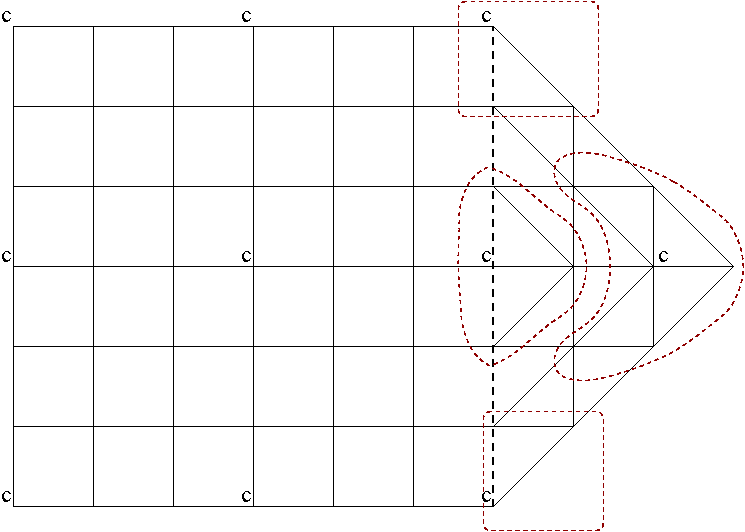}
  \caption{Structured square region attached to an unstructured triangular region.  The structure/unstructured interface is given by a dark dashed line. A \textit{c} denotes the location of a {\it Cpt}.  Red dashed lines encircle unstructured aggregates.  }
  \label{square+triangle}
\end{figure}

Figure~\ref{square+triangle} provides a simple illustration
of an unstructured triangular region attached to a $7 \times 7$ structured
region.  In Figure~\ref{square+triangle} a subset of vertices are labelled
with a `c' to denote a possible choice of coarse points denoted as {\it Cpts}.  The {\it Cpts}
set refers to a subset of fine mesh vertices that are chosen by a 
{\it classical} AMG algorithm to define the mesh vertices of the coarse mesh.
Notice that within structured regions, the {\it Cpts} have been defined in a
standard structured fashion. Ideally, it would be attractive to apply a 
standard AMG algorithm with no software modifications to coarsen
and define grid transfers for unstructured regions. However, the resulting
grid transfers stencils at co-located vertices must match their 
structured region counter-parts.  This means that the same set of
three {\it Cpts} should be chosen by the structured algorithm and the
unstructured algorithm along the interface in our Figure~\ref{square+triangle}
example and that the interpolation coefficients along the interface
be chosen in a very specific way.

In this paper, we do not employ classical AMG for unstructured
regions, but instead use the simpler plain aggregation variant of 
smoothed aggregation AMG method (SA)~\cite{Vanek1996a}. 
With both smoothed aggregation and plain aggregation multigrid, the coarsening
procedure is the same. In particular, 
coarsening is performed by aggregating together sets of fine vertices
as opposed to identifying {\it Cpts}. Each aggregate is essentially formed
by choosing a root vertex and including all of the root's neighbors
that have not already been included in another aggregate. Loosely, one can think of the aggregate root point as a {\it Cpt}. In Figure~\ref{square+triangle},
four aggregates in the unstructured region are depicted with dashed red lines.
To enforce the consistency of the {\it Cpts} choice at the interface,
the unstructured aggregation software must be changed so that it
initially chooses root points and aggregates associated with structured coarsening.
In our standard coarsening software, aggregation occurs in stages that are 
pipelined together. Each stage applies a specific algorithm that 
might only aggregate a subset of fine mesh vertices and then pass the 
partially-aggregated mesh to the next stage (that attempts to add more 
aggregates). Staging is a practical way to combine different aggregation 
algorithms with different objectives to ensure that all mesh vertices
are eventually aggregated. To accommodate structured/unstructured 
interfaces, a new aggregation stage was devised to start the aggregation 
process. This new stage only aggregates vertices on interfaces and chooses
root nodes in a structured fashion (employing a user-defined coarsening 
rate). Aggregates are chosen so that no interface vertices remain unaggregated
after this stage. Once this new stage completes, the standard unstructured
aggregation stages can proceed without further modification.
Notice that coarsening of structured and unstructured regions
can proceed fully in parallel (with no need for communication between the regions) as
processors responsible for unstructured regions redundantly coarsen/aggregate
the interface using the new devised aggregation stage while structured
regions also coarsen the interface using a standard structured coarsening
scheme. Since both structured and unstructured regions employ structured aggregation
along the mesh interface, matching {\it Cpts} are guaranteed.

Not only should coarsening be consistent along interfaces, but
interpolation coefficients at co-located vertices should match those 
produced by the structured regions.  For plain aggregation, multigrid this
will be the case as long as the structured region grid transfers
use the same methodology of piecewise constant basis functions.
Specifically, the corresponding plain aggregation interpolation basis functions 
are just piecewise constants for most applications. 
As the plain aggregation basis functions do not rely on the coefficients of 
the discretization matrix, each region's version of an interpolation 
stencil for a common interface will coincide exactly in the plane aggregation
case. This will not generally be true for more sophisticated AMG
schemes such as smoothed aggregation where the interpolation coefficients
depend on the discretization matrix coefficients. Effectively, a different
algorithm is used to generate the interpolation coefficients and so 
there is no reason why interpolation stencils should match those produced with linear
interpolation.
In this paper, we avoid this issue by only considering plain aggregation AMG
for unstructured regions in conjunction with piecewise constant interpolation (as opposed to linear
interpolation) for structured regions. However, we have identified two relatively 
straight-forward options both involving
some form of post-processing to the grid transfer operators. One possibility
is that a subset of processors 
communicate/coordinate with
each other to arrive at one common interpolation stencil for each unknown on a shared interface.
Obviously, this requires communication and is somewhat tedious to implement.
The second possibility is that linear basis functions always define interpolation
along interfaces between structured and unstructured regions. 
In this case, communication can be avoided by employing a
post-processing procedure within the unstructured grid transfer algorithm
to calculate (and overwrite) the appropriate interpolation operator along its 
interfaces.  
We omit the details but indicate that all the required information 
(coarse grid point locations and fine grid point locations) is already available
within our software framework.

To complete the discussion, we highlight some implementation
aspects associated with incorporating these pre- and post-processing changes
into a code such as {\muelu} which is based on a factory design, where
different classes must interact with different objects (e.g.,
aggregates, grid transfer matrices) needed to construct the multigrid
hierarchy. In particular, parameter lists are used to enter algorithm
choices and application specific data. In our context, the application
must indicate the following for each processor via parameter list entries:
\begin{itemize}
\item whether or not it owns a structured region or an unstructured region
\item the dimensions and coarsening rate for processors owning structured 
      regions
\item the dimensions and coarsening rate of each neighboring structured
      region for processors owning unstructured regions
\end{itemize}
Further, processors owning unstructured regions, that border structured
regions, must still provide structured region information for 
structured interfaces. This includes a  list of neighboring regions 
and the mapping of mesh nodes to regions as introduced in Figure~\ref{regionPerGIDWithGhosts}.

With the proper user-supplied information, {\muelu} assigns a
hybrid factory to address the prolongators. This hybrid factory 
includes an internal switch to then invoke either a structured region
grid transfer factory or an unstructured region grid transfer factory.
The hybrid factory essentially creates the grid transfer matrix object,
allowing the sub-factories to then populate this matrix object with 
suitable entries. It is this hybrid factory that invokes the 
aggregation process that starts with the interface aggregation stage
for unstructured regions. It is also responsible for the
post-processing (i.e., the updating of the prolongator matrix rows 
corresponding to interface rows) for the unstructured regions.
In this way, the standard structured factories and standard unstructured
factories require virtually no modifications, as these are mostly
confined to the hybrid factory. More information about {\muelu}'s factory
design can be found in \cite{BergerVergiat2019a}.

\section{Numerical Results}
Computational experiments are performed to highlight the equivalence
between MG cycles employing either composite operators or region 
operators as described by the Lemmas/Theorems presented earlier.
This is followed by experiments
to illustrate performance benefits of structured MG. Finally, we conclude this section
with an investigation demonstrating a structured region approach that also 
incorporates a few unstructured sub-domains. All the experiments that follow
can be reproduced using Trilinos at commit 86095f3d93e.


\subsection{Region MG Equivalence}
\label{sec:EquivalenceOfRegionMG}
To assess the equivalence of structured region MG to standard structured MG 
(without regions and region interfaces), we study a two-dimensional 
Laplace problem discretized with a $7$-point stencil on two different meshes,
a square $730\times730$ mesh and a rectangular $700\times720$ mesh.
The problem is run on~$9$ MPI ranks for the region solver and run in serial
for standard structured MG.
Here, we employ MG as a solver (not
as a preconditioner within a Krylov method), and the iteration is terminated 
when the relative residual drops below~$10^{-12}$.

The structured MG scheme employs a standard fully assembled matrix (i.e., a composite matrix
in this paper's terminology). It uses a coarsening rate of $3$ in each coordinate
direction and linear interpolation defines the grid transfer. The multigrid hierarchy consists of $4$ levels. 
Specifically, the hierarchy
mesh sizes from finest to coarsest for the square mesh are 
$730 \times 730$, $244 \times 244$, 
$82 \times 82$, and $28 \times 28$. Notice that all of these meshes correspond to $3^k +1$
points in each coordinate direction. Our software does not require these specific mesh
sizes, but this is needed to demonstrate exact equivalence. That is, both the composite
MG and the region MG must coarsen identically.  For the rectangular mesh, sizes
are not chosen so that the coarsening is identical (i.e., the number of 
vertices in each mesh dimension do not correspond to $3^k +1$). Thus, we expect
some small residual history differences for the rectangular mesh.
Fully structured multigrid is implemented
in Trilinos/{\muelu} using an option referred to as \emph{structured uncoupled aggregation}.
For the region MG hierarchy on the other hand, the mesh is partitioned into~$9~(= 3 \times 3)$ 
regions, where each region is assigned to one MPI rank. In this case, the square domain multigrid 
hierarchy for each processor's sub-mesh or region mesh is $244 \times 244$, $82 \times 82$, $28 \times 28$, and $9 \times 9$. In each coordinate direction, the overall finest mesh appears to have 
$732~(= 3$~processors $~\times~244 $ per processor) mesh points, which is not equal to
the $730$ mesh points used for the fully structured composite MG cycle. However, one must
keep in mind that 2 vertices are replicated along a mesh line in a coordinate direction (due to region the interfaces). Again,
these carefully chosen sizes are to enforce an identical coarsening procedure 
for the two MG solvers (and thus satisfy the conditions of the Lemmas/Theorems presented
earlier), as opposed to a hard requirement of the software. The region
multigrid method also uses a structured aggregation option to implement this type of 
structured coarsening. 

\Tabref{tab:EquivalenceOfResidualHistories} 
\begin{table}[h!]
\centering
\caption{Residual histories to study the equivalence of the structured region MG scheme to a classical structured MG}
\label{tab:EquivalenceOfResidualHistories}
\subfigure[$730\times730$ square mesh]{
\label{tab:EquivalenceSquareMesh}
\scriptsize
\begin{tabular}{c||cc||cc||cc}
~ & \multicolumn{2}{c||}{Jacobi} & \multicolumn{2}{c||}{\GS} & \multicolumn{2}{c}{\Cheby}\\
$\#$its. & Structured & $9$ Regions & Structured & $9$ Region & Structured & $9$ Regions\\
\hline
0  & 1.00000000e+00 & 1.00000000e+00 & 1.00000000e+00 & 1.00000000e+00 & 1.00000000e+00 & 1.00000000e+00\\
1  & 1.77885821e-02 & 1.77885821e-02 & 1.34144214e-02 & 1.34395087e-02 & 1.42870540e-02 & 1.42868592e-02\\
2  & 3.09066249e-03 & 3.09066249e-03 & 1.22727384e-03 & 1.23709339e-03 & 9.93752447e-04 & 9.93713870e-04\\
3  & 6.17432509e-04 & 6.17432509e-04 & 1.27481334e-04 & 1.29627870e-04 & 1.21921975e-04 & 1.21914771e-04\\
4  & 1.29973612e-04 & 1.29973612e-04 & 1.41133381e-05 & 1.45165400e-05 & 1.58413729e-05 & 1.58401012e-05\\
5  & 2.81812370e-05 & 2.81812370e-05 & 1.61878817e-06 & 1.69088891e-06 & 2.11105538e-06 & 2.11083642e-06\\
6  & 6.22574415e-06 & 6.22574415e-06 & 1.89847271e-07 & 2.02561731e-07 & 2.86037857e-07 & 2.86000509e-07\\
7  & 1.39312700e-06 & 1.39312700e-06 & 2.26276959e-08 & 2.48757453e-08 & 3.92564304e-08 & 3.92500462e-08\\
8  & 3.14666393e-07 & 3.14666393e-07 & 2.73250326e-09 & 3.13452182e-09 & 5.44989750e-09 & 5.44879379e-09\\
9  & 7.15836477e-08 & 7.15836477e-08 & 3.33798476e-10 & 4.06768456e-10 & 7.65555357e-10 & 7.65361045e-10\\
10 & 1.63770972e-08 & 1.63770972e-08 & 4.12201997e-11 & 5.46524944e-11 & 1.08974518e-10 & 1.08939546e-10\\
11 & 3.76413472e-09 & 3.76413472e-09 & 5.14512205e-12 & 7.64221900e-12 & 1.57581213e-11 & 1.57516868e-11\\
12 & 8.68493274e-10 & 8.68493274e-10 & 6.49387222e-13 & 1.11538919e-12 & 2.32197807e-12 & 2.32077246e-12\\
13 & 2.01044350e-10 & 2.01044350e-10 &                & 1.69735837e-13 & 3.49742848e-13 & 3.49514354e-13\\
14 & 4.66714466e-11 & 4.66714466e-11 &                &  &  & \\
15 & 1.08616953e-11 & 1.08616953e-11 &                &  &  & \\
16 & 2.53347464e-12 & 2.53347464e-12 &                &  &  & \\
17 & 5.92132868e-13 & 5.92132868e-13 &                &  &  & \\
\end{tabular}
}
\subfigure[$700\times720$ rectangular mesh]{
\label{tab:EquivalenceRectangularMesh}
\scriptsize
\begin{tabular}{c||cc||cc||cc}
~ & \multicolumn{2}{c||}{Jacobi} & \multicolumn{2}{c||}{\GS} & \multicolumn{2}{c}{\Cheby}\\
$\#$its. & Structured & $9$ Regions & Structured & $9$ Region & Structured & $9$ Regions\\
\hline
0  & 1.00000000e+00 & 1.00000000e+00 & 1.00000000e+00 & 1.00000000e+00 & 1.00000000e+00 & 1.00000000e+00\\
1  & 1.78374178e-02 & 1.77971728e-02 & 1.34028366e-02 & 1.34057178e-02 & 1.26092241e-02 & 1.25980465e-02\\
2  & 3.09747239e-03 & 3.08750444e-03 & 1.22692052e-03 & 1.22958855e-03 & 7.39937462e-04 & 7.40632616e-04\\
3  & 6.17958674e-04 & 6.15974350e-04 & 1.27486109e-04 & 1.28178073e-04 & 7.93385189e-05 & 7.96677401e-05\\
4  & 1.29899263e-04 & 1.29526261e-04 & 1.41232878e-05 & 1.42759476e-05 & 9.07488160e-06 & 9.15976761e-06\\
5  & 2.81258416e-05 & 2.80574257e-05 & 1.62135195e-06 & 1.65159920e-06 & 1.06848944e-06 & 1.08744092e-06\\
6  & 6.20516379e-06 & 6.19293768e-06 & 1.90317494e-07 & 1.95946605e-07 & 1.28512584e-07 & 1.32547397e-07\\
7  & 1.38672740e-06 & 1.38463243e-06 & 2.27023402e-08 & 2.37209815e-08 & 1.57501731e-08 & 1.65970557e-08\\
8  & 3.12830389e-07 & 3.12499063e-07 & 2.74346365e-09 & 2.92685712e-09 & 1.96757719e-09 & 2.14457022e-09\\
9  & 7.10802795e-08 & 7.10369390e-08 & 3.35333590e-10 & 3.68648440e-10 & 2.51098105e-10 & 2.87885059e-10\\
10 & 1.62430334e-08 & 1.62406131e-08 & 4.14287275e-11 & 4.75775741e-11 & 3.28456275e-11 & 4.04047421e-11\\
11 & 3.72913854e-09 & 3.73041913e-09 & 5.17289942e-12 & 6.32676763e-12 & 4.41956012e-12 & 5.94633832e-12\\
12 & 8.59490959e-10 & 8.60260047e-10 & 6.53051744e-13 & 8.72272622e-13 & 6.13278643e-13 & 9.15663966e-13\\
13 & 1.98754318e-10 & 1.99060540e-10 &  &  &  & \\
14 & 4.60939586e-11 & 4.62011981e-11 &  &  &  & \\
15 & 1.07170764e-11 & 1.07525542e-11 &  &  &  & \\
16 & 2.49746150e-12 & 2.50886608e-12 &  &  &  & \\
17 & 5.83206191e-13 & 5.86815903e-13 &  &  &  & \\
\end{tabular}
}
\end{table}
reports 
residual histories 
using Jacobi, {\GS}, and {\Cheby} as relaxation methods (1 pre- and 1 post-relaxation per level) in conjunction
with a direct solve on the coarsest level. In all cases, an identical right hand side and 
initial guess are used. 
Since the damped Jacobi smoother (which uses $\omega = .6$) only involves matrix-vector products and the true 
composite matrix diagonal, 
the residual histories match exactly for the square mesh. The square mesh residual histories are also
nearly identical with the {\Cheby} smoother, though there are small differences between the computed Chebyshev
eigenvalue intervals (whose calculation employs different random vectors).
In the case of the {\GS} relaxation, residual histories are still close, but do show slight differences.
This is due to the parallelization of {\GS}. As composite MG is run in serial,
it employs a true {\GS} algorithm while parallel region MG uses processor based (or domain decomposition based) 
{\GS}.  Specifically, applying {\GS} on a matrix row associated with a node in region~$\region{\dom}{i}$ on 
region interface~$\Boundss{i}{j}$ requires off-diagonal entries to represent the connections to neighboring nodes.
However, one (or more) neighboring nodes reside in the neighboring region~$\region{\dom}{j}$ and, thus, their matrix entries 
are not accessible for the {\GS} smoother.  The method does compute
the true composite residual before the {\GS} iteration, but only solution changes local to its
region are reflected in residual updates that occur within the smoother. 
Something similar occurs with composite MG {\GS} relaxation  in parallel, though the 
nature of its processor sub-domains are a bit different from those associated with regions. 
Even though the algorithms differ, one can see that the residual histories are close 
and only separate somewhat more significantly after more than 
$10$ orders of magnitude reduction in the residual.
The results for the rectangular mesh mirror those for the square mesh. The residual differences between
the standard composite MG and region MG are generally a tiny bit further from each other in this case
as the coarsening schemes for the two algorithms are no longer identical.

\subsection{Multigrid performance}

Region-based MG is motivated by potential performance gains when compared to a 
classical unstructured AMG method.  In the region-based case, one can exploit the regular 
structure of the mesh when designing both the data structure and implementing the key kernels 
used within the MG setup and {\vcycle} phases to avoid less indirect addressing and to reduce
the overall memory bandwidth requirements.

Our region MG is implemented in {\muelu}, which is part of the {Trilinos} framework. {Trilinos}
and {\muelu} have been designed and optimized for the type of fully unstructured meshes that 
might arise from a finite element discretization of a PDE problem. The underlying matrix data
structure is based on the Compressed Row Sparse format~\cite{Barret1994a} which can address these
types of general sparse unstructured data. At present, our region MG software is in its initial 
stages 
and so it utilizes these same underlying unstructured data formats for matrices and vectors. 
Thus, it has not been optimized for structured grids. Interestingly, we are able to 
demonstrate some performance gains in the case of PDE systems, even with the current software 
limitations. We begin first with some Poisson results 
and then follow this with elasticity experiments where 
significant gains are observed. In both cases, linear finite elements with hexahedral
elements are used to construct the linear systems.

For both the Poisson and the elasticity experiments, the problem setup is as follows.
Each region performs coarsening by a rate of $3$, until three levels have been formed.
On the coarsest region-level~$\numLevelsRegion-1$, we then apply AMG as a coarse level solver 
as outlined in \secref{sec:CoarseLevelSolver}.  Depending on the problem size on the finest 
level, $1-3$ rounds of additional coarsening will be performed algebraically until the coarse operator
of the AMG hierarchy has less than $900$ rows and can be tackled by a direct solver.
On all levels~$\indLevel\in\{0,1,\hdots,\numLevels-2\}$, but the coarsest, damped Jacobi
smoothing is employed using a damping parameter of $.67$. That is, both the region hierarchy
and the {\it coarse-solver} AMG hierarchy use the same smoother settings.
On the coarsest region-level~$\numLevelsRegion-1$, each MPI rank only owns a few rows,
so a repartitioning/rebalancing step is performed before constructing the AMG coarse level solver to avoid
having a poorly balanced AMG coarse solve that requires a significant amount of communication.


To avoid confusion, we now use the term pure AMG to describe the standard AMG approach 
(without any levels using a region format) that is used for the comparisons. 
The pure AMG hierarchy uses the 
same smoother settings employed
for the region multigrid method
as well as the 
same total number of levels~$\numLevels$ 
(counting both the region/structured levels and {\it coarse-solver} AMG levels).
As with region MG, a direct solver is applied on the coarsest level.
In all cases where AMG is employed, level transfer operators 
are constructed using SA-AMG \cite{Vanek1996a} with {\muelu}'s uncoupled aggregation 
and a prolongator smoothing damping parameter~$\omega = 4/3$.
To counteract poor load balancing during coarsening,
we repartition such that each MPI rank at least owns $800$ rows and that the relative mismatch in size between all subdomains is less than $10\%$.
Partitioning is perform via multi-jagged coordinate partitioning using Trilinos' {\zoltanTwo} package\footnote{\url{https://trilinos.github.io/zoltan2.html}}.
Since our examples focus on a direct comparison of region MG and AMG,
we apply the MG scheme as a solver without any outer Krylov method. 
Of course, application codes will often invoke MG as a preconditioner within a Krylov method. We report timings for the both the MG hierarchy setup and for the solution phase 
of the algorithm.

Table~\ref{tab:PoissonTimingsRegionMGvsAMG} and Table~\ref{tab:ElasticityTimingsRegionMGvsAMG_Jac}
present the timings.
\begin{table}
\centering
\caption{Region MG {\vs} AMG for three-dimensional Poisson example: configuration and performance}
\label{tab:PoissonTimingsRegionMGvsAMG}
\begin{tabular}{c|c|c|ccc|ccc}
Mesh & $\nproc$ & $\numLevels$ & \multicolumn{3}{c|}{Structured MG} & \multicolumn{3}{c}{Pure Algebraic MG}\\
    nodes & ~ & $\numLevelsRegion/\numLevelsAMG$ ($\numLevels$) & $\#$its & Setup & {\vcycle} & $\#$its & Setup & {\vcycle}\\
\hline
$82^3$  & $27$    & 3/2 (3) & $13$ & $\unit{0.0728}{\second}$ & $\unit{0.193}{\second}$ & $13$ & $\unit{0.117}{\second}$ & $\unit{0.242}{\second}$\\
$163^3$ & $216$   & 3/2 (4) & $13$ & $\unit{0.104 }{\second}$ & $\unit{0.241}{\second}$ & $13$ & $\unit{0.176}{\second}$ & $\unit{0.273}{\second}$\\
$325^3$ & $1728$  & 3/3 (5) & $13$ & $\unit{0.352 }{\second}$ & $\unit{0.428}{\second}$ & $13$ & $\unit{0.581}{\second}$ & $\unit{0.400}{\second}$\\
$622^3$ & $12167$ & 3/3 (6) & $13$ & $\unit{0.386 }{\second}$ & $\unit{0.425}{\second}$ & $13$ & $\unit{0.711}{\second}$ & $\unit{0.423}{\second}$
\end{tabular}
\end{table}
\begin{table}
\centering
\caption{Region MG {\vs} AMG for three-dimensional elasticity example: configuration and performance for Jacobi smoother}
\label{tab:ElasticityTimingsRegionMGvsAMG_Jac}
\begin{tabular}{c|c|c|ccc|ccc}
Mesh & $\nproc$ & $\#$levels & \multicolumn{3}{c|}{Structured MG} & \multicolumn{3}{c}{Pure Algebraic MG}\\
nodes & ~ & $\numLevelsRegion/\numLevelsAMG$ ($\numLevels$) & $\#$its & Setup & {\vcycle} & $\#$its & Setup & {\vcycle}\\
\hline
$82^3$  & $27$    & 3/2 (4) & $22$ & $\unit{0.333}{\second}$ & $\unit{1.94}{\second}$ & $35$ & $\unit{2.46}{\second}$ & $\unit{4.23}{\second}$\\
$163^3$ & $216$   & 3/3 (5) & $21$ & $\unit{0.423}{\second}$ & $\unit{1.97}{\second}$ & $33$ & $\unit{2.78}{\second}$ & $\unit{4.34}{\second}$\\
$325^3$ & $1728$  & 3/3 (5) & $21$ & $\unit{0.697}{\second}$ & $\unit{2.38}{\second}$ & $32$ & $\unit{3.54}{\second}$ & $\unit{4.92}{\second}$\\
$622^3$ & $12167$ & 3/4 (6) & $20$ & $\unit{1.199}{\second}$ & $\unit{2.63}{\second}$ & $32$ & $\unit{3.92}{\second}$ & $\unit{5.06}{\second}$
\end{tabular}
\end{table}
These tests were performed in parallel on Cori\footnote{\url{https://docs.nersc.gov/systems/cori/}}
at the National Energy Research Scientific Computing Center (NERSC), Berkeley, CA. 
The mesh sizes as well as parallel resources are given in the first two 
columns of each table. 
The column entitled ``mesh nodes'' denotes the number of grid nodes in the cube-type mesh.
The number of MPI ranks~$\nproc$ is increased at the same rate as the mesh size, yielding a weak scaling type of experiment.
For the region MG algorithm, the number of MPI ranks also denotes the number of regions,
such that the number of unknowns per region is kept constant across all experiments at $\approx 20k$ unknowns per MPI rank.

The gains for the Poisson problem correspond to about a factor of two in the setup phase. It is important
to recall that many of the key computational kernels (e.g., the matrix-matrix multiply) employ the same code for the region MG 
and for pure AMG. These setup gains come primarily from a faster process to generate grid transfers and having somewhat fewer nonzeros
within the coarse level matrices. Without doubt, the most time consuming kernel on larger core counts comes from repartitioning
the matrix supplied to the coarse AMG solver. This repartitioning reduces communication costs in constructing the coarse
AMG hierarchy, but it comes with a high price.  While the actual data transfer associated with rebalancing requires some communication, 
the great bulk of this repartitioning time involves the cost associated with using Trilinos' framework to set up the communication 
data structure (which includes some neighbor discovery process).  It is important to notice that when solving a sequence of linear systems
on the same mesh (e.g., within a nonlinear solution scheme or within a time stepping algorithm), this communication data structure 
remains the same throughout the sequence
\footnote{This would not necessarily be true for an AMG scheme that uses a strength-of-connection method that effectively alters the matrix-graph based on the matrix's nonzero values.}.
Thus, it should be possible to form this data structure just once and reuse it over the entire sequence,
drastically reducing this communication cost.

The elasticity 
results exhibit more than a factor of three improvement in the setup phase and a factor of two in the solve phase, even without using kernels geared
toward structured grids. In the case of AMG setup, this is mostly due to the lower number of coarse operators nonzeros.
This is reflected in multigrid operator complexities (which measures the ratio of the total number of nonzeros in the discretization
matrices on all levels versus the number of nonzeros in the finest level matrix. In the region case it is under 1.1 (which includes nonzeros
associated with {\it coarse-solver} AMG levels). In the pure AMG case it is over 1.4. 
Additionally, there are some savings in that no communication is required while constructing the region part of the hierarchy, though 
once again there are costs associated with the coarse AMG setup.
For the solve phase, the benefits come from having less nonzeros and also requiring fewer iterations, which is 
due to the fact that linear interpolation is the better grid transfer than that provided by SA-AMG for this problem.

\subsection{Multigrid kernel performance}

While the current structured region code is unoptimized, we have started experimenting with alternative multigrid kernels
outside of the Trilinos package. 
%
In this section we illustrate the potential gains that may be possible even while retaining a matrix data structure best
suited for fully unstructured grids. Specifically, timing comparisons are made between the multigrid matrix-matrix multiply kernel 
from our standard unstructured AMG package, Trilinos/{\muelu}, and a special purpose one written for two dimensional structured meshes.
This special purpose matrix-matrix multiply 
also 
requires a small amount of additional information (e.g., number of grid points in the coordinate directions for each region).
In all cases, the kernels produce the same results (with the exception of slight numerical rounding variations).
The only difference is that the new kernel leverages the structured grid layout.
While one might consider designing new data structures to support structured kernels, we are currently evaluating
tradeoffs. 
%
Using the same unstructured data structures greatly facilitates the integration and maintenance of the new structured capabilities within our predominantly unstructured AMG package,
though it may somewhat curb or limit the performance gains attained by the structured kernels.

For the matrix-matrix multiplication 
the underlying matrix data structure consists of two integer arrays and one double precision array associated with the 
compressed row matrix format~\cite{Barret1994a}.
One of the integer arrays consists of pointers to the starting location (within the other two arrays) of the data corresponding to a matrix row.
The other two arrays hold column indices and matrix values for the nonzeros.
While all three arrays are still passed to the 
matrix-multiply kernel,
one nice benefit of the structured algorithms is that access to the two integer arrays can be limited.
In particular, all the data within the integer arrays can be inferred or deduced once the structured stencil pattern and grid layout are known.
This ultimately reduces memory access and allows for a number of other optimizations. See \cite{Bergen2004a} for some examples.

To demonstrate the matrix-multiply gains, we evaluate the matrix triple product or Galerkin projection step within the multigrid setup phase
corresponding to 
$$
\bar{A} = R A P .
$$
A two dimensional mesh is considered along with a perfect factor of three coarsening in each coordinate direction. 
For the unstructured {\muelu} implementation, the product $A P$ is first formed using a two-matrix multiplication procedure.
The product of $R$ and the result of the first two-matrix multiplication is then performed to arrive at the desired result.
For the structured implementation, the triple product is formed directly. That is, explicit formulas have been determined
(using a combination of Matlab, Mathematica, and pre/post processing programs)
for each of $\bar{A}$'s entries.
Specifically, there are four sets of formulas for rows of $\bar{A}$ corresponding to each of the four mesh corners.
There are an additional four sets of formulas for the four mesh sides (excluding the corners).
Finally, there is one last set of formulas for the mesh interior. As noted above, the integer arrays are not used in the evaluation of these formulas.

\begin{figure}
  \centering
\includegraphics[trim = 0.3in 0.3in 0.3in 0.3in, clip, height = 5.1cm,width = 5.1cm]{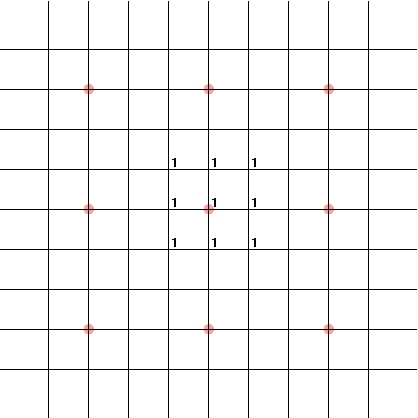}
~~~~~
\includegraphics[trim = 0.3in 0.3in 0.3in 0.3in, clip, height = 5.1cm,width = 5.1cm]{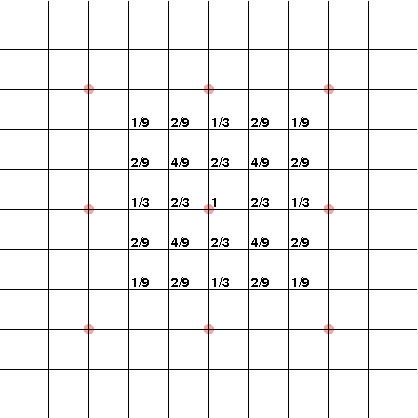}
  \caption{One grid transfer column stencil associated with the central coarse point using piecewise constants (left) and linear interpolation (right). Only a portion of the mesh is shown and circles denote coarse mesh points.}
  \label{rap stencils}
\end{figure}

Three different structured functions have been developed. One corresponds to the use of piecewise constant grid transfers;
another is for geometric grid transfers on a regular uniform mesh; the third allows for general grid transfers
(which have the same sparsity pattern as the geometric grid transfers but allow for general coefficient values).
An interior basis function stencil (or column) is depicted in Figure~\ref{rap stencils} for the piecewise constant case and for the ideal geometric case.
In these two contexts, the coefficients of $R$, and $P$ do not need to be accessed as they are known ahead of time and have been included in the explicit formulas.
In the general situation, the double precision arrays for $R$ and $P$ must be accessed to perform the triple product.
In all cases, $A$ is assumed to have a nine point stencil within the interior.
Stencils along the boundary have the same structure where entries are dropped if they are associated with points that extend outside of the mesh.

Table~\ref{rap timing} illustrates some representative serial timings.  The reported mesh sizes refer to the coarse mesh. The 
corresponding fine mesh is given by $(3 n_x - 2) \times (3 n_y -2)$ for a coarse $n_x \times n_y$ mesh.
Here, one can see that the structured versions are generally an order of magnitude faster than the unstructured Trilinos/{\muelu} kernel.
\begin{table}
\begin{center}
\caption{Timings (in seconds) for different triple-matrix product kernels.
\textsc{9 pt Basis} (\textsc{25 pt Basis}) indicates $9$ ($25$) point basis functions for $P$ and $R$.
\textsc{Const}, \textsc{Geo}, and \textsc{Generic} denote the structured triple product for piecewise constant,
ideal geometric, and general grid transfers.
} \label{rap timing}
\begin{tabular}{c|cc||ccc}
coarse          & \multicolumn{2}{c||}{\textsc{9 pt Basis}} & \multicolumn{3}{c}{\textsc{25 pt Basis}} \\
mesh size       & {\muelu} & \textsc{Const} & {\muelu} & \textsc{Generic} & \textsc{Geo} \\ \hline
$140 \times 36 $&  .0024 &   .0001        &  .0109 &     .0012         &    .0009    \\
$140 \times 180$&  .0124 &   .0006        &  .0572 &     .0061         &    .0046     \\
$700 \times 180$&  .0726 &   .0070        &  .2944 &     .0320        &     .0240     \\
$700 \times 900$&  .3702 &   .0356        & 1.4786 &     .1606        &     .1208     
\end{tabular}
\end{center}
\end{table}
These timings correspond to the core multiply time (excluding a modest amount of time needed in Trilinos
to pre/post process data to pre-compute additional information needed for parallel computations). 
As no inter-region communication is required (due to Theorem~\ref{rap theorem}), the structured serial run 
times are representative of parallel run times when one region is assigned to each processor. 
Given the fact the triple product is one of the most costly AMG setup kernels and the fact that
the Trilinos matrix-matrix multiply has been optimized many times over the years, these $10x$ gains
are significant. 

It should be noted, however, that we have not integrated the improved triple products into our
framework. In particular, we have not yet developed efficient 3D formulas, which is somewhat labor intensive
to perform properly. Additionally, we still have several framework decisions concerning how different
structured grid cases are addressed and merged within our generally unstructured AMG package.

\subsection{Multigrid for hybrid structured/unstructured meshes}
\begin{figure}
\begin{center}
    \begin{tikzpicture}[yslant=-0.05]
\draw (0,0,3)--(3,0,3)--(3,3,3)--(0,3,3)--cycle;
\draw (3,0,0)--(3,3,0)--(3,3,3)--(3,0,3)--cycle;
\draw (3,3,0)--(3,3,3)--(0,3,3)--(0,3,0)--cycle;
\foreach \x in {1,2}
{
\draw(3,0,\x)--(3,3,\x); 
\draw(3,\x,0)--(3,\x,3); 
\draw(0,\x,3)--(3,\x,3); 
\draw(\x,0,3)--(\x,3,3); 
\draw(\x,3,0)--(\x,3,3); 
\draw(0,3,\x)--(3,3,\x); 
}
\node at (0.5,0.5,3) {0};
\node at (1.5,0.5,3) {1};
\node at (2.5,0.5,3) {2};
\node[rotate=0, xslant=0.0, yslant=0.95] at (3,0.5,2.5) {\scalebox{0.6}[1.0]{2}};%
\node[rotate=0, xslant=0.0, yslant=0.95] at (3,0.5,1.5) {\scalebox{0.6}[1.0]{5}};%
\node[rotate=0, xslant=0.0, yslant=0.95] at (3,0.5,0.5) {\scalebox{0.6}[1.0]{8}};%
\node at (0.5,1.5,3) {9};
\node at (1.5,1.5,3) {10};
\node at (2.5,1.5,3) {11};
\node[rotate=0, xslant=0.0, yslant=0.95] at (3,1.5,2.5) {\scalebox{0.6}[1.0]{11}};%
\node[rotate=0, xslant=0.0, yslant=0.95] at (3,1.5,1.5) {\scalebox{0.6}[1.0]{14}};%
\node[rotate=0, xslant=0.0, yslant=0.95] at (3,1.5,0.5) {\scalebox{0.6}[1.0]{17}};%
\node at (0.5,2.5,3) {18};
\node at (1.5,2.5,3) {19};
\node at (2.5,2.5,3) {20};
\node[rotate=0, xslant=0.0, yslant=0.95] at (3,2.5,2.5) {\scalebox{0.6}[1.0]{20}};%
\node[rotate=0, xslant=0.0, yslant=0.95] at (3,2.5,1.5) {\scalebox{0.6}[1.0]{23}};%
\node[rotate=0, xslant=0.0, yslant=0.95] at (3,2.5,0.5) {\scalebox{0.6}[1.0]{26}};%
\node[rotate=-3, xslant=0.95, yslant=0.00] at (0.5,3,2.5) {\scalebox{1.2}[0.5]{18}};
\node[rotate=-3, xslant=0.95, yslant=0.00] at (1.5,3,2.5) {\scalebox{1.2}[0.5]{19}};
\node[rotate=-3, xslant=0.95, yslant=0.00] at (2.5,3,2.5) {\scalebox{1.2}[0.5]{20}};
\node[rotate=-3, xslant=0.95, yslant=0.00] at (0.5,3,1.5) {\scalebox{1.2}[0.5]{21}};
\node[rotate=-3, xslant=0.95, yslant=0.00] at (1.5,3,1.5) {\scalebox{1.2}[0.5]{22}};
\node[rotate=-3, xslant=0.95, yslant=0.00] at (2.5,3,1.5) {\scalebox{1.2}[0.5]{23}};
\node[rotate=-3, xslant=0.95, yslant=0.00] at (0.5,3,0.5) {\scalebox{1.2}[0.5]{24}};
\node[rotate=-3, xslant=0.95, yslant=0.00] at (1.5,3,0.5) {\scalebox{1.2}[0.5]{25}};
\node[rotate=-3, xslant=0.95, yslant=0.00] at (2.5,3,0.5) {\scalebox{1.2}[0.5]{26}};
\end{tikzpicture}
    \scalebox{1}[1]{\includegraphics[scale=0.21, trim=0 90pt 0 0]{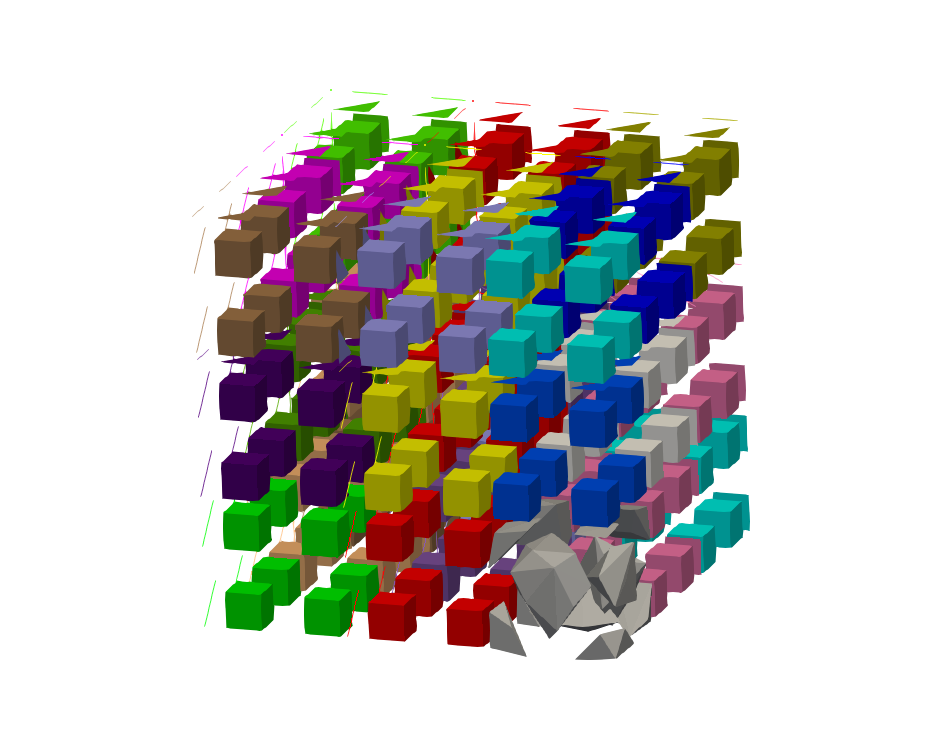}}
\end{center}
\caption{On the left, a visualization of a $3\times3\times3$ Region layout on a cube. On the right, an example of the region aggregates, with region 2 unstructured.}
\label{fig:27ProcLayout}
\end{figure}
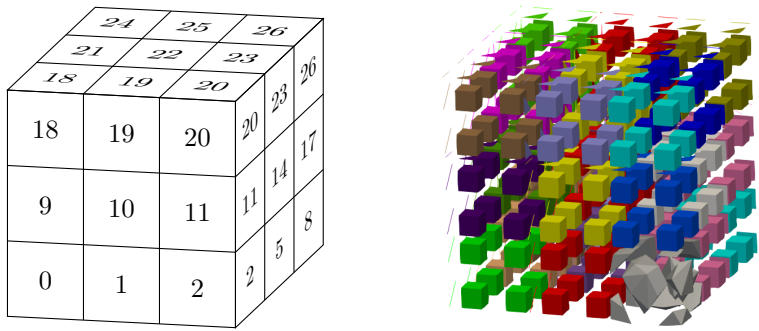

To demonstrate the flexibility of the proposed region MG scheme to handle semi-structured meshes containing unstructured regions
we consider a $3\times3\times3$ region setup with different regions flagged as either structured or unstructured.
The region layout is illustrated in Figure~\ref{fig:27ProcLayout}
along with a visualization of the aggregates when one region, region $2$, is treated as unstructured.
For the numerical tests, we solve a 3D Poisson equation with a $7$-point stencil on a
$100\times100\times100$ mesh cube using a 3-level {\wcycle} and piecewise constant interpolation
for both the structured multigrid and for the unstructured region AMG.
Presently, our implementation only properly addresses a structured/unstructured region combination
using piecewise constant interpolation (i.e., the Lemmas presented in this paper are satisfied).
Proper extensions for linear interpolation (discussed in Section~\ref{sec:StructuredUnstructured})
are planned for a a refactored version of the software.
Two iterations of Symmetric {\GS} are used as the pre and post smoothers, and the coarse grid is solved with a direct solve.
The problem is solved to a tolerance of $10^{-6}$.
Table~\ref{structuredunstructured} shows iteration counts when different regions are marked as unstructured,
and the remaining regions are structured.

We see that the introduction of unstructured regions does have a small impact on the convergence rate of the method,
with more unstructured regions resulting in slightly more iterations, up to the limit of all regions being treated as unstructured.
This is likely a result of suboptimal aggregates being formed along the interfaces due to the forced matching of aggregates between neighboring regions.
We have observed that this effect is more pronounced when the coarsening rate in the structured regions differs from the coarsening rate of the unstructured region (in experiments not shown in this paper).
Here, the structured regions used a coarsening rate of 3 and the unstructured regions have an approximate coarsening rate of 3 as well.


\begin{table}
\begin{center}
\caption{
Iteration counts for various structured/unstructured setups. The regions are setup in a $3\times3\times3$ format.
For structured/unstructured testing, we solve a 3D Laplace equation on a $100\times100\times100$ cube.
Two iterations of Symmetric {\GS} are used as the pre smooth and post smooth for a
3-level W-cycle multigrid iteration with piecewise constant interpolation.
} \label{structuredunstructured}
\begin{tabular}{|c|c|} \hline
    Region Layout                  & Iterations \\ \hline
    AMG with no region formatting  & 17 \\ \hline
    no unstructured regions        & 15 \\ \hline
    no structured regions          & 18 \\ \hline
    Front Face unstructured        & 17 \\ \hline
    Back Face unstructured         & 17 \\ \hline
    Top Face unstructured          & 17 \\ \hline
    Bottom Face unstructured       & 16 \\ \hline
    Left Face unstructured         & 17 \\ \hline
    Right Face unstructured        & 16 \\ \hline
    Eight Corners unstructured     & 16 \\ \hline
    Region 2 unstructured          & 15 \\ \hline
    Region 13 unstructured         & 16 \\ \hline
    Region 24 unstructured         & 15 \\ \hline
    Regions 2, 13, 24 unstructured & 16 \\ \hline
\end{tabular}
\end{center}
\end{table}

\section{Concluding remarks}
\label{sec:Conclusion}
We have presented a generalization of the HHG idea to a semi-structured framework.
Within this framework, the original computational domain is decomposed into regions
that only overlap at inter-region interfaces.  Unknowns along region interfaces are
replicated so that each region has its own copy of the solution along its interfaces.
This facilitates the use of structured grid kernels within a multigrid algorithm when
regions are structured.  We have presented a mathematical framework to represent
this region decomposition. The framework allows us to precisely define components
of  a region multigrid algorithm and understand the conditions by which such a region
multigrid algorithm is identical to a traditional multigrid algorithm. Using this framework,
we illustrate how a region multigrid hierarchy can be constructed without requiring inter-region
communication in some cases. We have also presented some ideas towards making the use
of such a region multigrid solver less invasive for application developers. These ideas
exploit transformations that define conversions between a region representation and a more
traditional representation for vectors and matrices. We also illustrated how such a
multigrid solver can account for some unstructured regions within the domain.
Finally, we have presented some
evidence of the potential of such an approach in terms of computational performance.

\bibliographystyle{abbrv}
\bibliography{bib_region_multigrid.bib}

\begin{thebibliography}{10}

\bibitem{Barret1994a}
R.~Barret, M.~Berry, T.~F. Chan, J.~Demmel, J.~Donato, J.~Dongarra,
  V.~Eijkhout, R.~Pozo, C.~Romine, and H.~v.~d. Vorst.
\newblock {\em {Templates for the Solution of Linear Systems: Building Blocks
  for Iterative Methods}}.
\newblock SIAM, Philadelphia, PA, USA, 1994.

\bibitem{Bergen2006a}
B.~K. Bergen, T.~Gradl, F.~H{\"u}lsemann, and U.~R{\"u}de.
\newblock {A Massively Parallel Multigrid Method for Finite Elements}.
\newblock {\em Computing in Science \& Engineering}, 8(6):56--62, 2006.

\bibitem{Bergen2004a}
B.~K. Bergen and F.~H{\"u}lsemann.
\newblock Hierarchical hybrid grids: data structures and core algorithms for
  multigrid.
\newblock {\em Numerical Linear Algebra with Applications}, 11(2-3):279--291,
  2004.

\bibitem{Bergen2007a}
B.~K. Bergen, G.~Wellein, F.~H{\"u}lsemann, and U.~R{\"u}de.
\newblock {Hierarchical hybrid grids: achieving TERAFLOP performance on large
  scale finite element simulations}.
\newblock {\em International Journal of Parallel, Emergent and Distributed
  Systems}, 22(4):311--329, 2007.

\bibitem{MueLuURL}
L.~Berger-Vergiat, C.~A. Glusa, J.~J. Hu, M.~Mayr, P.~Ohm, A.~Prokopenko, C.~M.
  Siefert, R.~S. Tuminaro, and T.~A. Wiesner.
\newblock {The MueLu Multigrid Framework}.
\newblock \url{https://trilinos.github.io/muelu.html}, 2020.

\bibitem{BergerVergiat2019a}
L.~Berger-Vergiat, C.~A. Glusa, J.~J. Hu, M.~Mayr, A.~Prokopenko, C.~M.
  Siefert, R.~S. Tuminaro, and T.~A. Wiesner.
\newblock {MueLu User's Guide}.
\newblock Technical Report SAND2019-0537, Sandia National Laboratories,
  Albuquerque, NM (USA) 87185, 2019.

\bibitem{Briggs2000a}
W.~L. Briggs, V.~E. Henson, and S.~F. McCormick.
\newblock {\em {A Multigrid Tutorial}}.
\newblock SIAM, 2nd edition, 2000.

\bibitem{Dendy2010a}
J.~E. Dendy and J.~D. Moulton.
\newblock Black box multigrid with coarsening by a factor of three.
\newblock {\em Numerical Linear Algebra with Applications}, 17(2-3):577--598,
  2010.

\bibitem{Dubey20143217}
A.~Dubey, A.~Almgren, J.~Bell, M.~Berzins, S.~Brandt, G.~Bryan, P.~Colella,
  D.~Graves, M.~Lijewski, F.~L{\"o}ffler, B.~O'Shea, E.~Schnetter, B.~V.
  Straalen, and K.~Weide.
\newblock A survey of high level frameworks in block-structured adaptive mesh
  refinement packages.
\newblock {\em J. of Par. and Distr. Comput.}, 74(12):3217 -- 3227, 2014.

\bibitem{hypre}
R.~Falgout, J.~Jones, and U.~Yang.
\newblock The design and implementation of hypre, a library of parallel high
  performance preconditioners.
\newblock In A.~Bruaset and A.~Tveito, editors, {\em Numerical Solution of
  Partial Differential Equations on Parallel Computers}, volume~51 of {\em
  Lecture Notes in Computational Science and Engineering}. Springer, Berlin,
  2006.

\bibitem{Gmeiner2013a}
B.~Gmeiner, T.~Gradl, F.~Gaspar, and U.~R{\"u}de.
\newblock {Optimization of the multigrid-convergence rate on semi-structured
  meshes by local Fourier analysis}.
\newblock {\em Computers \& Mathematics with Applications}, 65(4):694--711,
  2013.

\bibitem{Gmeiner2015b}
B.~Gmeiner, M.~Huber, L.~John, U.~R{\"u}de, and B.~I. Wohlmuth.
\newblock {A quantitative performance study for Stokes solvers at the extreme
  scale}.
\newblock {\em Journal of Computational Science}, 17(3):509--521, 2016.

\bibitem{Gmeiner2012a}
B.~Gmeiner, M.~Mohr, and U.~R{\"u}de.
\newblock {Hierarchical Hybrid Grids for Mantle Convection: A First Study}.
\newblock In {\em 2012 11th International Symposium on Parallel and Distributed
  Computing}, pages 309--314, 2012.

\bibitem{Gmeiner2015c}
B.~Gmeiner, U.~R{\"u}de, H.~Stengel, C.~Waluga, and B.~I. Wohlmuth.
\newblock {Performance and Scalability of Hierarchical Hybrid Multigrid Solvers
  for Stokes Systems}.
\newblock {\em SIAM Journal on Scientific Computing}, 37(2):C143--C168, 2015.

\bibitem{Hackbusch1994a}
W.~Hackbusch.
\newblock {\em {Iterative Solution of Large Sparse Systems of Equations}},
  volume~95 of {\em Applied Mathematical Sciences}.
\newblock Springer, 1994.

\bibitem{Henshaw20087469}
W.~Henshaw and D.~Schwendeman.
\newblock Parallel computation of three-dimensional flows using overlapping
  grids with adaptive mesh refinement.
\newblock {\em J. of Comp. Phys.}, 227(16):7469 -- 7502, 2008.

\bibitem{Lee03asynchronousfast}
B.~Lee, S.~Mccormick, B.~Philip, and D.~Quinlan.
\newblock Asynchronous fast adaptive composite-grid methods: Numerical results.
\newblock {\em SIAM J. Sci. Comput.}, 25:2003, 2003.

\bibitem{Philip20122277}
B.~Philip and T.~Chartier.
\newblock Adaptive algebraic smoothers.
\newblock {\em J. of Comp. and Appl. Math.}, 236(9):2277 -- 2297, 2012.

\bibitem{Prokopenko2016b}
A.~Prokopenko, C.~M. Siefert, J.~J. Hu, M.~Hoemmen, and A.~Klinvex.
\newblock {Ifpack2 User's Guide 1.0}.
\newblock Technical Report SAND2016-5338, Sandia National Laboratories, 2016.

\bibitem{Saad2003a}
Y.~Saad.
\newblock {\em {Iterative Methods for Sparse Linear Systems}}.
\newblock SIAM, Philadelphia, PA, USA, 2003.

\bibitem{octrees}
R.~Sampath and G.~Biros.
\newblock A parallel geometric multigrid method for finite elements on octree
  meshes.
\newblock {\em SIAM J. Sci. Comput.}, 32(3):1361--1392, 2010.

\bibitem{uintah}
J.~Schmidt, M.~Berzins, J.~Thornock, T.~Saad, and J.~Sutherland.
\newblock Large scale parallel solution of incompressible flow problems using
  {Uintah} and {Hypre}.
\newblock In {\em Cluster, Cloud and Grid Computing (CCGrid), 2013 13th
  IEEE/ACM International Symposium on}, pages 458--465, May 2013.

\bibitem{TrOoSc00}
U.~Trottenberg, C.~W. Oosterlee, and A.~Schuller.
\newblock {\em Multigrid}.
\newblock Academic Press, 2000.

\bibitem{Vanek1996a}
P.~Van\v{e}k, J.~Mandel, and M.~Brezina.
\newblock {Algebraic Multigrid By Smoothed Aggregation For Second And Fourth
  Order Elliptic Problems}.
\newblock {\em Computing}, 56:179--196, 1996.

\end{thebibliography}

\end{document}